\documentclass[a4paper,fleqn]{cas-dc}

\usepackage[authoryear]{natbib}
\usepackage[subrefformat=parens]{subcaption}

\usepackage{lineno}
\newcommand*\patchAmsMathEnvironmentForLineno[1]{%
  \expandafter\let\csname old#1\expandafter\endcsname\csname #1\endcsname
  \expandafter\let\csname oldend#1\expandafter\endcsname\csname end#1\endcsname
  \renewenvironment{#1}%
     {\linenomath\csname old#1\endcsname}%
     {\csname oldend#1\endcsname\endlinenomath}}% 
\newcommand*\patchBothAmsMathEnvironmentsForLineno[1]{%
  \patchAmsMathEnvironmentForLineno{#1}%
  \patchAmsMathEnvironmentForLineno{#1*}}%
\AtBeginDocument{%
\patchBothAmsMathEnvironmentsForLineno{equation}%
\patchBothAmsMathEnvironmentsForLineno{align}%
\patchBothAmsMathEnvironmentsForLineno{flalign}%
\patchBothAmsMathEnvironmentsForLineno{alignat}%
\patchBothAmsMathEnvironmentsForLineno{gather}%
\patchBothAmsMathEnvironmentsForLineno{multline}%
}
\usepackage{ulem}
\newcommand{\Add}[1]{\textcolor{black}{#1}}	%
\newcommand{\Erase}[1]{\if0{#1}\fi}	%

\newcommand{\AddTwo}[1]{\textcolor{black}{#1}}	%
\newcommand{\EraseTwo}[1]{\if0{#1}\fi}	%

\usepackage{amsmath}   % for math commands
\usepackage{mathtools} % for math commands
\usepackage{cleveref}  % for clever reference

\ExplSyntaxOn
\keys_set:nn { stm / mktitle } { nologo }
\ExplSyntaxOff

\crefname{equation}{Eq.}{Eqs.}% 
\crefname{figure}{Fig.}{Figs.}%

\def\tsc#1{\csdef{#1}{\textsc{\lowercase{#1}}\xspace}}
\tsc{WGM}
\tsc{QE}
\tsc{EP}
\tsc{PMS}
\tsc{BEC}
\tsc{DE}
\begin{document}
\let\WriteBookmarks\relax
\def\floatpagepagefraction{1}
\def\textpagefraction{.001}
\shorttitle{Optimization on Planning of Trajectory and Control of Autonomous Berthing and Unberthing for the Realistic Port Geometry}
\shortauthors{Miyauchi et~al.}

\title [mode = title]{Optimization on Planning of Trajectory and Control of Autonomous Berthing and Unberthing for the Realistic Port Geometry}                      
% \tnotemark[1,2]

% \tnotetext[1]{This document is the results of the research
%   project funded by the National Science Foundation.}

% \tnotetext[2]{The second title footnote which is a longer text matter
%   to fill through the whole text width and overflow into
%   another line in the footnotes area of the first page.}

\let\printorcid\relax % Remove ORCID footnote

\author[1]{Yoshiki Miyauchi}
\cormark[1]
% \fnmark[1]
% \ead{cvr_1@tug.org.in}
% \ead[url]{www.cvr.cc, cvr@sayahna.org}

\credit{Methodology, Software, Writing - original draft, visualization}

\address[1]{Department of Naval Architecture and Ocean Engineering, Graduate School of Engineering, Osaka University, 2-1 Yamadaoka, Suita, Osaka 565-0971, Japan}
\ead{yoshiki_miyauchi@naoe.eng.osaka-u.ac.jp}

\author[1,2]{Ryohei Sawada}
\credit{Methodology, Writing - review \& editing}
\address[2]{National Maritime Research Institute, 6-38-1, Shinkawa, Mitaka, Tokyo 181-0004, Japan}

\author[3,4]{Youhei Akimoto}
\credit{Software, Writing - review \& editing}
\address[3]{Faculty of Engineering, Information and Systems, University of Tsukuba, 1-1-1 Tennodai, Tsukuba, Ibaraki 305-8573, Japan}
\address[4]{RIKEN Center for Advanced Intelligence Project, 1-4-1 Nihonbashi, Chuo-ku, Tokyo 103-0027, Japan}
\author[1]{Naoya Umeda}
\credit{Supervision, Funding acquisition}

\author[1]{Atsuo Maki}
\credit{Conceptualization, Writing - review \& editing}
\cormark[1]
\ead{maki@naoe.eng.osaka-u.ac.jp}

\cortext[cor1]{Corresponding author}

\begin{abstract}
To realize autonomous shipping, autonomous berthing and unberthing are some of the technical challenges. In the past, numerous research have been done on the optimization of trajectory planning of berthing problems. However, these studies assumed only a simple berth and did not consider obstacles. Optimization of trajectory planning on berthing and unberthing in actual ports must consider the spatial constraints and maintain sufficient distance to obstacles. The main contributions of this study are as follows: (i) a collision avoidance algorithm based on the ship domain which has variable size by the ship speed is proposed, to include the spatial constraints to optimization; (ii) the effect of wind disturbance is taken into account to the trajectory planning to make a feasible trajectory based on the capacity limit of actuators; (iii) showing that the optimization method for berthing is also eligible for the unberthing, which has been almost neglected; (iv) waypoints are included to the optimization process, to make optimization easier on practical applications. The authors tested the proposed method on two existing ports.  The proposed method performed well on both the berthing and the unberthing problem and optimized the control input and the trajectory while avoiding collision with the complex obstacles.
\end{abstract}

% \begin{graphicalabstract}
% \includegraphics{figs/grabs.pdf}
% \end{graphicalabstract}

% \begin{highlights}
% \item \Add{A collision avoidance algorithm based on the ship domain was proposed.}
% \item \Add{The optimization of berthing control was performed with the proposed method.}
% \item \Add{The proposed method can obtain a collision-free optimal control input at real port.}
% \item \Add{The waypoint navigation was included in the optimization to make it more robust.}
% \end{highlights}


\maketitle

\sloppy

\section{Introduction}\label{sec:intro}
Realizing autonomous ship navigation, autonomous berthing is one of the major technical challenges.  
Operations such as berthing and unberthing at a narrow fairway are the most stressful for the navigator. This is because \textit{``The difficulty is not only due to the increase in the number of maneuvering measures, but also because the ship is forced to maneuver in a severe surrounding environment, including spatial constraints such as breakwaters and berths, reduced maneuverability due to low-speed navigation, and increased effects of wind and currents.''} \citep{Seta2004} Thus,  berthing and unberthing are relatively difficult than that of navigation of open-sea. Thus automation of these kinds of maneuvering is useful not only from the viewpoint of reducing the operational cost but also from reducing the workload on the navigator.

In solving the online control of the autonomous berthing, it is effective to use the predefined trajectory and control input as a reference for tracking control obtained by the optimization problem. The process to generate the predefined trajectory is called trajectory planning. In this paper, we defined that trajectory as time series of state of the ship, including geometric paths, with some or all dynamic constraints of the ship, by referring to the definition of \citet{Vagale2021}: \textit{``path planning is typically defined within purely geometric space, whereas trajectory planning, or trajectory generation, involves geometric paths endowed with temporal properties, e.g., to incorporate dynamics.'' }  Trajectory planning is generally difficult to directly apply to the control algorithm due to its computation requirements, which is longer than the computation time required for real-time control\Erase{caused by the non-linearity of the maneuvering motion}. However, when applying real-time control such as PID control or a black box approach such as neural networks, an optimal trajectory can be used as a reference trajectory to configure a feasible control that is optimized in terms of time and other indicators while satisfying constraints based on actuator capacity limitations and maneuverability.

\Add{Another application of trajectory planning is the evaluation tool of ship's actuator's capability design on berthing. If the trajectory planning tool can sufficiently incorporate the ship's dynamics, actuator capability, spatial constraints, and external disturbance, the designer of a ship can evaluate the ship's berthing capability on specific port geometry and wind condition by utilizing the trajectory planning tool. This would be valuable on the design stage of both autonomous and manned ships.}

%離着桟の経路最適化の先行研究をかく
\subsection{Related Research}\label{related}
In the past, numerous research have been done on the optimization of trajectory planning and control of berthing problems \citep{Hasegawa1993, Ahmed2015,Maki2020,Maki2020b}. However, most of these studies have assumed only a simple straight berth and do not consider surrounding obstacles. To perform berthing and unberthing in actual ports, it is necessary to avoid static obstacles such as bridge piers, berths other than the destination and anchored ships, as well as moving obstacles such as passing ships.

Only a few research have included the complex geometric constraints of the port to trajectory planning of berthing. Those research have used a two-stage method which combined the path planning using graph search and optimal control problem solver for trajectory generation: optimization of trajectory planning with polygonal constrains which combined the graph-search and convex optimization \citep{Martinsen2019, Martinsen2020}; combined the hybrid A$^{*}$ search algorithm as an initial guess to an optimal control problem solver under the external disturbance \citep{Bitar2020}; the combination of lattice-based motion planning and receding horizon improvement\citep{Bergman2020}. They created an initial trajectory using the lattice-based motion planning which transformed the optimal control problem to classical graph search problem. As a second stage, they conducted optimization based improvement of the initial trajectory using receding horizon improvement while the ship maintaining a fixed safe distance to obstacles.

% These research with complex constrains are valuable, however the differences with present paper are are follows.
\Add{The difference between present work and those similar previous research \citep{Martinsen2020,Bitar2020,Bergman2020} are as follows. First, They }have adopted the two-stage method, which \Erase{discrete the continuous state} \Add{translate the problem} to graph search \Add{problem} \Erase{algorithm which obtains} \Add{to obtain an sub-optimal} initial guess \Add{at the first stage. The Second stage}  optimize the initial guess with an optimization problem solver. \Add{The proposed method is a single-stage method that does not require an initial guess. Second, they used a simplified dynamic model for the subject ship's maneuvering system: azimuth thruster-like thrust vectoring model. The present work assumed a more common configuration: rudder and propeller system with side thrusters. To incorporate this common configuration, more complex maneuvering models are required than previous research. Finally,} these research \Erase{also} have only considered the constant \Erase{radius of} safety region to the obstacle\Add{, whereas the proposed method used dynamic safety region. }

Trajectory planning on other kinds of restricted waters have been reported on: global trajectory planning and collision avoidance of an autonomous surface vessel at narrow ferry passage by safety region and collision region around moving and static obstacles with applying the hybrid dynamic window method \citep{Serigstad2018}; Navigation in restricted channels using reinforcement learning was done by \citet{Amendola2020}.

% 船以外の経路最適化の研究をかく

On the optimization of berthing and unberthing trajectory, it is important to consider whether the trajectory avoids the obstacle and whether the trajectory is separated from the obstacle by a sufficient distance for safety. The problem of how much distance a ship should keep from obstacles, was investigated by \citet{fujii1971} in the 1970s as \textit{effective domain}; and followed by \citet{goodwin1975}, which defined \textit{ship domain} as \textit{``the effective area around a ship which a navigator would like to keep free with respect to other ships and stationary objects.''} Since then, various models have been proposed and applied to the collision avoidance method. \citet{Szlapczynski2017} provides a comprehensive review of these ship domain models. However,  many models did not address the geographical constraints because the ship domain was mainly used to represents the safe distance between ships. Effective domain at restricted waters was proposed by \citet{Yamanouchi1972} based on observation of navigation outside the port of Yokohama, Japan, and showed that in restricted waters, the ship domain is smaller than open seas, which about three-fourths in the longitudinal direction and about half in the width direction.  Ship domain based on traffic at strait proposed by \citet{Wang2016} that is proportional to the ship length and speed. The coefficients are determined based on the Vessel Traffic Information System (VTIS) database of ships transiting the Singapore Strait by using a genetic algorithm. \citet{Pietrzykowski2008} proposed a ship fuzzy domain with the shape of the narrow fairway as a parameter. \citet{Hansen2013} proposed the ship domain on restricted waters based on AIS data measured on the Fehmarnbelt between Denmark and Germany and \citet{Jensen2013} followed that work by estimating the traffic separating scheme's efficiency on the bridge pier of Fehmarnbelt. However, in our literature survey, the ship domain designed to apply to berthing or unberthing problems cannot be found. \Add{To approach the berth wall at the terminal phase of berthing, distance to the obstacle will be much smaller than the ship domain stated above. A novel ship domain applicable for berthing is required.}

\subsection{Objective and Scope}
The objective of this study is to include several major factors necessary for applying optimization of trajectory planning on both berthing and unberthing to practical application, such as spatial constraints and external disturbance. Trajectory planning on unberthing has not been focused on previous research compare to berthing; however, it is equally important to realize the fully autonomous ship. The main contributions of this study are as follows: (i) a new collision avoidance algorithm based on the ship domain which has variable size by the ship speed was proposed, to include the spatial constraints to \Add{single-stage} optimization of trajectory planning \Add{ with complex maneuvering model}; (ii) effect of wind disturbance was taken into account to the trajectory planning to make a feasible trajectory based on the capacity limit of actuators; (iii) showing that the optimization method for berthing is also eligible for the unberthing, which has been almost neglected on the field on trajectory planning; (iv) waypoints are included to the optimization process, to make optimization easier on practical applications. Major environmental conditions related to berthing and unberthing consist of wind, wave, and current. However, disturbance caused by wave and current were neglected in this study, assuming that the vessel enters a port sheltered by breakwaters. 

By adding the idea of ship domain to the collision avoidance algorithm, the proposed method enables the search for a trajectory that maintains an appropriate distance from the obstacle. In this study, the obstacle was limited to a static obstacle, but the algorithm can be extended to a dynamic obstacle.

This study extends the study of \citet{Maki2020} the Co-variance Matrix Adaption Evolution Strategy (CMA-ES) for the optimization process. \Add{The study of \citet{Maki2020} only assumed a simple trajectory planning on straight line berth, hence the method had several limitations on practical use such as: can not handle multiple, arbitrary shaped berth; safety distance to obstacles was not considered; external disturbance was not considered. Those limitations are addressed in this study. By incorporating both spatial constraints and wind disturbance, the proposed method can use for both reference trajectory generator of autonomous berthing and evaluation tool of ship's design on berthing capability.} 

The subject ship is a large single-rudder, single-propeller ship with side thrusters, which can berth and unberth without assistance by tug boats. The proposed method requires only the initial and desired end state, control constraints, and spatial constraints of obstacle geometry as computational conditions. To validate the proposed algorithm, we also performed calculations on multiple ports.

The rest of the paper is organized as follows: Section 2 describes the optimization method of trajectory planning and mathematical modeling of maneuver. The optimization method contains the collision avoidance algorithm, which is newly proposed in this study; Section 3 shows the test result on trajectory planning using a newly introduced optimization scheme with a collision avoidance algorithm. Section 3 also describes the geometry of ports which served as test cases; In Section 4, a discussion on the results shown in the previous section and future works are shown; Finally, section 5 concludes the study. 

\section{Method}\label{sec:method}
\begin{figure}
    \centering
    \includegraphics[width=\linewidth]{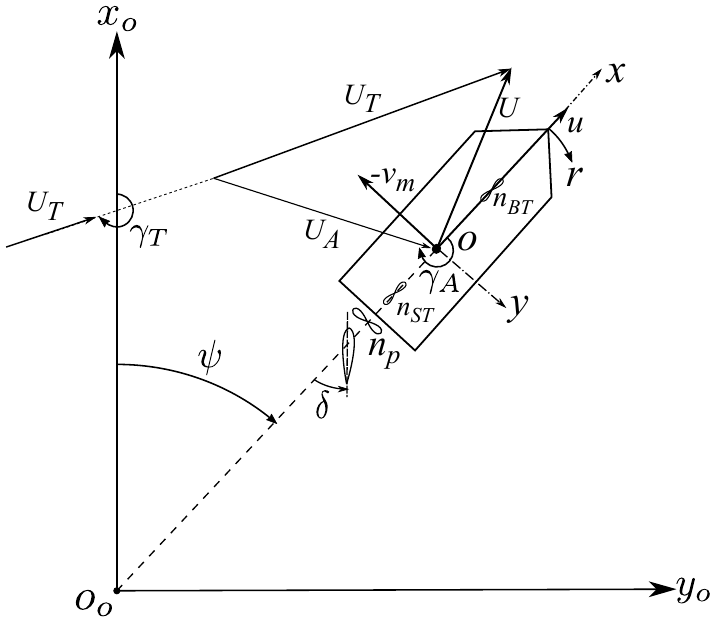}
    \caption{Coordinate System. \Add{The coordinate system consist by space-fixed system $o_{0}-x_{0}y_{0}$ and ship-fixed system $o-xy$.} \Add{Notations represent: $u,v_{m}$ are the velocity of $o-xy$ system;  $\delta,~n_{\mathrm{p}},~n_{\mathrm{BT}},~n_{\mathrm{ST}}$ are the rudder angle, the revolution of propeller, the bow thruster and the stern thruster; $\gamma_{T},~U_{T}$ are true wind direction and speed; $\gamma_{A},~U_{A}$ are apparent wind direction and speed.}}
    \label{fig:coordinate}
\end{figure}

In this study, berthing and unberthing are modeled as 3 degrees-of-freedom problem. The coordinate systems are space-fixed system $o_{0}-x_{0}y_{0}$ and ship-fixed system $o-xy$ which has its origin on midship.  The origin of $o_{0}-x_{0}y_{0}$ was set as the position of midship at designated berthing point or the start point of unberthing, true-north direction as the positive direction of $x_{0}$, east direction as the positive direction of $y_{0}$. \Cref{fig:coordinate} shows the coordinate systems in this study. State vector is $\mbox{\boldmath $x$}=(x_{0}, u, y_{0}, v_{m}, \psi, r)^\mathsf{T}  \in \mathbb{R}^{6}$, where $u,~v_{m} $ are the velocity of $o-xy$ system. The vector of control input is $\mbox{\boldmath $u$}=(\delta, n_{p}, n_{BT}, n_{ST})^\mathsf{T} $. $\delta,~n_{\mathrm{p}},~ n_{\mathrm{BT}},~n_{\mathrm{ST}}$ represent the rudder angle, the revolution of propeller, the bow thruster and the stern thruster, respectively. The positive direction of the rudder angle was set to the direction in which the heading angle 
$\psi$ increase by steering, positive rotation of propeller was set to the advancing direction, positive rotation of side thrusters was set to the direction which the induced force work to move the ship to the positive direction of ship-fixed coordinate $y$. 

Wind disturbance was considered as external force $\boldsymbol{\omega}$, which consist by the true wind direction and true wind speed $\boldsymbol{\omega}=(\gamma_{T},~U_{T})^\mathsf{T} $. However, due to the ship's maneuver, an apparent wind affects as the actual force acting on the hull. Hence the apparent wind was computed inside the mathematical model of maneuvering. The zero direction of true wind direction $\gamma_{T}$ was set to the direction in which the wind blows from the positive direction to the negative direction of $x_{0}$.

\subsection{Optimization of Trajectory planning}
\subsubsection{Objective Function}

In this study, the berthing and unberthing problems were modeled as time-minimizing problem. \Add{Aiming points of the objective function for berthing and unberthing problem are 1) to be collision-free, 2) to satisfy the desired end state, and 3) to minimize the time. Indeed, 1) and 2) are essential constraints that cannot be compromised, but in this study, they are not constraints but included in the objective function. Candidates for evaluation indices other than collision-free and satisfaction of end state can be: minimization of time; minimization of distance; minimization of energy consumption; and ride quality improvement. Regarding those candidates, minimizing time is easy to express in quantitative terms, and minimizing time leads to a short and efficient trajectory. On the other hand, minimizing the distance itself has little practical meaning than minimizing time. Minimizing energy is also not much meaningful since berthing and unberthing account for only a small portion of a ship's energy consumption. The optimization of ride quality is useful for passenger ships, but we avoided it because of the lack of proper quantitative expression of the ride quality.} 

The state equation shown is \cref{eq:dynamics}:
    \begin{equation}
        \dot{\boldsymbol{x}}(t) = \boldsymbol{f}\big( \boldsymbol{x}(t), \boldsymbol{u}(t), \boldsymbol{\omega}(t) \big).
        \label{eq:dynamics}
    \end{equation}
On the berthing and unberthing problem, the initial and end state are given.  By expanding the objective function $J$ proposed by \citet{Maki2020}, $J$ in this study is:
% \begin{align}
%     \begin{aligned}
%      \label{eq:objfunc}
%     J=\sum_{i}^{6} \max \left[w_{\mathrm{pen}} \boldsymbol{w}_{\mathrm{dim}, i} \cdot\left\{\boldsymbol{x}_{\mathrm{des}, i}-\boldsymbol{x}_{i}(t_{\mathrm{f}})\right\}^{2},\right. \\
%     \left. \boldsymbol{w}_{\mathrm{dim}, i} \cdot \boldsymbol{x}_{\mathrm{tol}, i}^{2}\right] \cdot t_{\mathrm{f}} +w_{c} C  
%     \end{aligned}
% \end{align}
\begin{multline} \label{eq:objfunc}
    J =  w_c \cdot C 
    + t_\mathrm{f} \cdot \sum_{i=1}^{6} w_{\mathrm{dim},i}\bigg(x_{\mathrm{tol},i}^2 \mathbf{1}_{\{\lvert x_{\mathrm{des}, i}-x_{i}(t_{\mathrm{f}})\rvert \leq x_{\mathrm{tol},i} \}} \\ 
    + w_\mathrm{pen} \big(x_{\mathrm{des}, i}-x_{i}(t_{\mathrm{f}})\big)^{2} \mathbf{1}_{\{\lvert x_{\mathrm{des}, i}-x_{i}(t_{\mathrm{f}})\rvert > x_{\mathrm{tol},i} \}} \bigg)\enspace,
\end{multline}
where the difference $x_{\mathrm{des}, 5}-x_{5}(t_{\mathrm{f}})$ in angle should be treated as a value in $(-\pi, \pi]$ by taking into account the minimum difference angle and $\mathbf{1}_{\{A\}}$ is the indicator function that is $1$ if $A$ is true and $0$ otherwise. Here, subscript $i=1,2, \cdots, 6$ means a related term to the $i$-th component of the state vector.
The first term of $J$ is the penalty term for the collision with obstacles, and the second term is for satisfying the desired end state condition.  $\boldsymbol{x}_{\mathrm{int}}$ and $\boldsymbol{x}_{\mathrm{des}}$ are the vectors of given initial state and the desired end state.
 $\boldsymbol{x}_{\mathrm{tol}}$ is the tolerance vector to the difference between the actual and desired end state, which was introduced to prevent $t_{\mathrm{f}}$ from increase due to the pursuit with excessive precision of the end state during the optimization process, same as the previous study \citep{Maki2020}. 

In the case of unberthing, $u$ does not need to be in the certain range of tolerance, instead it should be better just above the certain exit speed. Thus, for the unberthing case, the condition of indicator function of \cref{eq:objfunc} replaced by  $u(t_{\mathrm{f}})\geq \boldsymbol{x}_{\mathrm{tol}, 2}$. $w_{\mathrm{pen}}$ is the penalty coefficient when the actual end state not satisfy the $\boldsymbol{x}_{\mathrm{tol}}$, $\boldsymbol{w}_{\mathrm{dim}}$ is the collection vector for state vector components which have different dimensions. On the second term, $C$ is the collision penalty function and $w_{\mathrm{c}}$ is the weight coefficient for collision penalty. Details of collision penalty function is shown on \cref{sec:collision}. List of coefficient and $\boldsymbol{x}_{\mathrm{tol}}$ are shown on  Table\ref{tab:list_wcoeffs}.

There is some discussion on how to configure $\boldsymbol{x}_{\mathrm{tol}}$. On the berthing speed in Japanese ports, \citet{Murakami2015} proposed a regression equation based on measurement. However, their study was conducted from the perspective of allowable values in port design, not from how much the ship needs to slow down. On the heading angle, \citet{Roubos2017} analyzed 555 records of berthing of tankers, bulkers, and container ships in the port of Rotterdam and found that the approach angle to the berth at the moment of impact is less than $1.5^{\circ}$. The tolerance $\boldsymbol{x}_{\mathrm{tol}}$ was set based on the tolerance values of the state at the time of berthing in these previous studies. The velocity component of $\boldsymbol{x}_{\mathrm{tol}}$ was calculated from the regression equation for berthing velocity $V$ which covers 90\% of the measured data  \citep{Murakami2015}: 
\begin{equation}
    \label{eq:murakami}
    V=0.279 \cdot GT^{-0.114},
\end{equation}
where $GT$ is the gross tonnage and was assumed as a ferry of approximately 10,000 tons, which operated on the subject port shown on section~\ref{sec:port_geo}.
The tolerance of $\psi$ was set to $\boldsymbol{x}_{\mathrm{tol},5}=1.0^{\circ}$, which is slightly smaller than the maximum value in \citet{Roubos2017}. The position's tolerance was set to $1.0\mathrm{m}$, which is small enough for the ship's length.

During the optimization, the optimization method searches for a combination of the time series $\boldsymbol{u}(t)$ and $t_{\mathrm{f}}$ that minimizes \cref{eq:objfunc}. The time series of the control input is kept constant during the period of $t_{c}=90$ s to prevent frequent input changes and to keep the number of variables to be searched finite. Hence time series of control inputs is discretized into $m$ segments. Thus the variable $\boldsymbol{X}$ to be optimized is follows:
% \begin{align}
%     \label{eq:controlvector}
%     \begin{aligned}
%     \boldsymbol{X} \equiv &\left(  t_{\mathrm{f}},  \delta_{1}, \cdots, \delta_{m}, n_{p,1}, \cdots, n_{p,m}, n_{BT, 1}, \cdots, n_{BT, m}, \right. \\
%     &\left. n_{ST, 1}, \cdots, n_{ST, m} \right)^\mathsf{T}  \in \mathbf{R}^{4 m+1}   
%     \end{aligned}
% \end{align}
\begin{multline}
    \label{eq:controlvector}
    \boldsymbol{X} \equiv \big(  t_{\mathrm{f}},  \delta_{1}, \cdots, \delta_{m}, n_{p,1}, \cdots, n_{p,m}, n_{BT, 1}, \cdots, n_{BT, m}, \\
    n_{ST, 1}, \cdots, n_{ST, m} \big)^\mathsf{T}  \in \mathbf{R}^{4 m+1}  \enspace.
\end{multline}
In this study $m=25$, hence the maximum time of simulation is $2250$ seconds.

Finally, the optimization problem can be expressed as follows:
\begin{align}
    \begin{aligned}
        \boldsymbol{X}_{\mathrm{opt}} &=\text {argmin} \ 
        J(t_{f}, \boldsymbol{u})\\
        \text { subject to : } x(0) &=x_{\text{int}} \\
        t_{\mathrm{f, min}}\leq & t_{\mathrm{f}} \leq t_{\mathrm{f, max}} \\
        |\delta(t)| & \leq \delta_{\max } \\
        |n_{\mathrm{p}}(t)| & \leq n_{\mathrm{p}\max } \\
        |n_{\mathrm{BT}}(t)| & \leq n_{\mathrm{BT}\max } \\
        |n_{\mathrm{ST}}(t)| & \leq n_{\mathrm{ST}\max } \enspace.
    \end{aligned}
\end{align}
The optimal trajectory is derived by substituting the optimal solution ${X}_{\mathrm{opt}}$ to the mathematical model of maneuver. The minimum and maximum value on the exploration domain of $\boldsymbol{X}$ are shown in \Cref{tab:minmax}. The limits were set to: the maximum rudder angle of common rudder system; the propeller revolution of advancing at 8 kn; a revolution which can induce the equivalent force of the lateral wind force at 30 m/s when using both of the bow and the stern thruster. 
\begin{table*}
    \centering
    \caption{list of weight coefficients and tolerance vectors of end state}
    \begin{tabular}{ccc}
        \hline Parameter & \multicolumn{2}{c} { Value } \\
        \hline \hline Condition & Berthing & Unberthing \\
        $\boldsymbol{x}_{\mathrm{tol}}$ & $(1.0 \mathrm{m}, \ 0.1 \mathrm{m/s}, \ 1.0 \mathrm{m}, \ 0.1 \mathrm{m/s}, \ 1.0 ^{\circ}, \ 0.0764^{\circ}/\mathrm{s})^\mathsf{T} $ & $(1.0 \mathrm{m},\  0.0 \mathrm{m/s}, \ 1.0 \mathrm{m}, \ 0.1 \mathrm{m/s}, \ 1.0 ^{\circ}, \ 0.764^{\circ}/\mathrm{s})^\mathsf{T}  $\\
        $w_{\mathrm{pen}}$ &  \multicolumn{2}{c}{$1.0\times10^{4}$}  \\
        $\boldsymbol{w}_{\mathrm{dim}}$ & \multicolumn{2}{c}{$\left(1 / w_{L}^{2}, \ 1 / w_{U}^{2}, \ 1 / w_{L}^{2}, \ 1 / w_{U}^{2}, \ \pi^{2}, \ w_{L}^{2} / w_{U}^{2}\right)^\mathsf{T} $} \\
        $w_{\mathrm{C}}$ & \multicolumn{2}{c}{$1.0\times10^{10}$}  \\
        $w_{L}$ & \multicolumn{2}{c}{$0.1L_{\mathrm{pp}}$} \\
        $w_{U}$ & $\boldsymbol{x}_{\mathrm{init}, 2}/2$ & $\boldsymbol{x}_{\mathrm{des}, 2}/2$\\
        \hline
    \end{tabular}
    \label{tab:list_wcoeffs}
\end{table*}

\begin{table}
    \centering
    \caption{Minimum and maximum of end time $t_{\mathrm{f}}$ and control inputs}
    \begin{tabular}{cc}
        \hline Parameter & Value \\
        \hline\hline $t_{\mathrm{f, min}}$, $\mathrm{(s)}$  & $630$ \\
        $t_{\mathrm{f, max}}$, $\mathrm{(s)}$   & $2250$ \\
        $|\delta_{\max}|$, $\mathrm{(degree)}$ & $35$ \\ 
        $|n_{\mathrm{p}, \max}|$,$\mathrm{(rps)}$ & $2.08$ \\ 
        $|n_{\mathrm{BT}, \max}$|,$\mathrm{(rps)}$ & $4.24$ \\ 
        $|n_{\mathrm{ST}, \max}$|,$\mathrm{(rps)}$ & $4.24$ \\ 
        \hline
    \end{tabular}
    \label{tab:minmax}
\end{table}

\subsubsection{Optimization Method: CMA-ES}
\begin{figure}
    \centering
    \includegraphics[width=\linewidth]{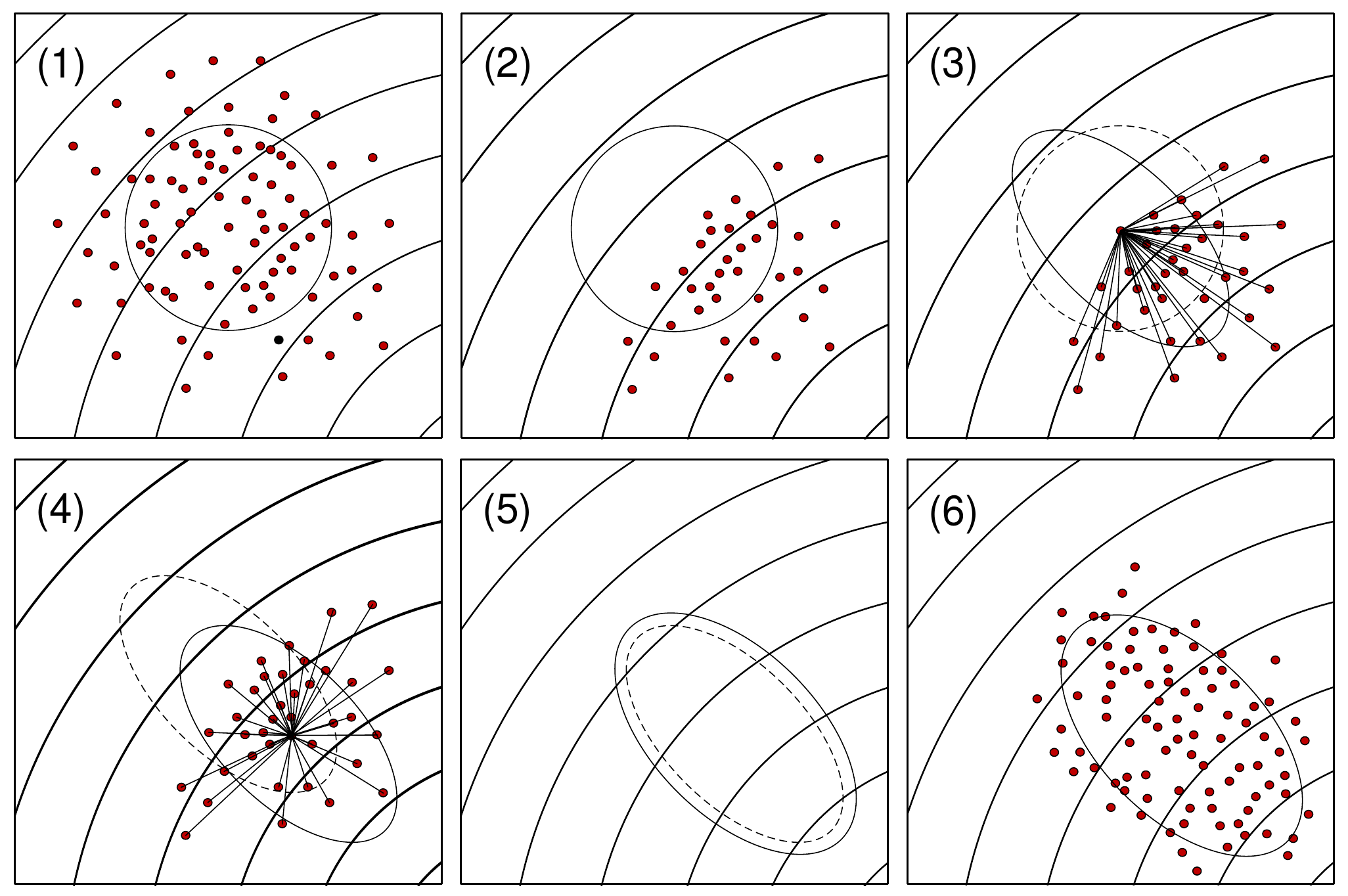}
    \caption{Schematic presentation of the CMA-ES procedure including (1) generating multiple candidate solutions, (2) evaluating and ranking the solutions based on the objective function, (3) updating the covariance matrix, (4) shifting the center of the distribution to a weighted mean vector, (5) updating the step size and (6) generating multiple candidates in the next step. This figure duplicates Fig.~2 in the literature \citep{Maki2020b}}
    \label{fig:cma}
\end{figure}

As in the previous study \citep{Maki2020, Maki2020b}, co-variance matrix adaption evolution strategy (CMA-ES) \citep{Hansen2006} adapted to box constrains \citep{Sakamoto2017} and with restart strategy \citep{Auger2005} was applied as optimization method. Fig.\ref{fig:cma} shows the procedure of optimization using CMA-ES. In previous research \citep{Maki2020}, CMA-ES showed the capability to obtain a feasible solution on optimal berthing problem as a time-minimizing problem while taking into account the collision penalty to berth with simple geometry. 

\begin{figure}
    \centering
    \includegraphics[width=\columnwidth]{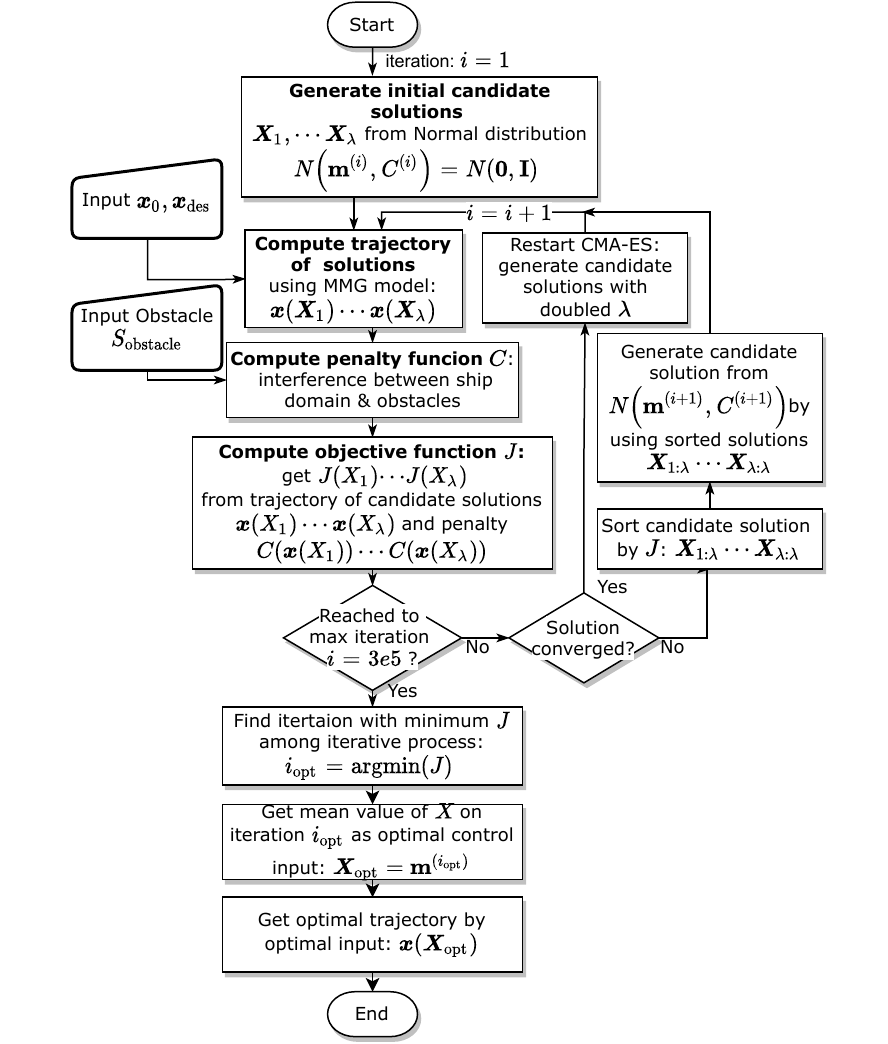}
    \caption{Flowchart of optimization computation using CMA-ES.}
    \label{fig:flowchart}
\end{figure}

\Add{The flow of optimization computation are as follows. \Cref{fig:flowchart} illustrates the flow. First, CMA-ES generated $\lambda$ candidate solutions $\boldsymbol{X}_{j=1 \cdots \lambda}$ from normal distribution $N(\mathbf{m}^{(i)}, \mathbf{C}^{(i)})$. Second, Generated candidate solutions were evaluated by objective function $J$ (\Cref{eq:objfunc}). 
To compute $J$, time history of state vector $\boldsymbol{x}(t)$ on solution $\boldsymbol{X_{j=1\cdots\lambda}}$ were computed by numerical simulation using MMG model (\Cref{sec:MMG}). Once the $\boldsymbol{x}$(t) for $\boldsymbol{X}_{j=1\cdots\lambda}$ was obtained, the penalty function $C$ was computed to represent interference between ship domain and obstacles for each solution. Details of $C$ are shown on \Cref{sec:collision}. 
Third, the $\mathbf{m}^{(i)}$ and $\mathbf{C^{(i)}}$ were updated using weighted average of candidate solutions. Updated parameters are used to generate candidate solutions of next iteration.
Candidate solutions were sorted in order of $J$ before taking weighted average, to generate better candidate in the next iteration. 
If the candidate solutions were converged during optimization, The CMA-ES was restated. On the detail of restart strategy, please see \citep{Maki2020b}}. In this study, the initial population size of CMA-ES is 20, and the max size is 240, while the population size was doubled when the restart occurred.
\Add{Those solution generation-evaluation-update iteration was continued until user defined maximum iteration number was satisfied. The max. iteration number was set to $3\times10^{5}$. And finally, when after the iteration was stopped, optimal solution $\boldsymbol{X}_{\text{opt}}$ was searched from result of iterative process. Again, the optimal trajectory $\boldsymbol{x}(\boldsymbol{X}_{\text{opt}})$ was derived by feeding $\boldsymbol{X}_{\text{opt}}$ to maneuvering simulation.}

\subsubsection{Collision avoidance Algorithm}\label{sec:collision}
In the previous research \citep{Maki2020}, collision avoidance was than by adding penalty function $C$ to objective function instead of imposing constraint to system.The penalty function $C$ was defined as the integration of instantaneous penetration length of the vertex of the rectangle surrounding the ship:
\begin{equation}
    C=\sum_{i=1}^{4} \int_{P_{i} \in C_{\text {berth }}}\left|Y_{i}-Y_{\text {berth }}\right| dt \enspace,
\end{equation}
where $P_{i}$ is the four vertices of the surrounding rectangle, $C_{\mathrm{berth}}$ is the area inside the berth, $Y_{i}$ denotes each  $y_{0}$ ordinate of $P_{i}$ and  $Y_{\text {berth }}$ represents the $y_{0}$ coordinate of the berth edge.
This collision avoidance method has several limitations for practical use: \Add{(a)} it can deal with only one berth at a time; \Add{(b)} the edge of a berth must be parallel to the $y_{0}$ direction; \Add{(c)} safety distance to obstacles, such as berths, breakwaters, buoys or anchoring vessels, is not considered. 

\Add{The limitation (c): lack of safety distance is the most critical limitation from the practical application perspective. A trajectory which does not maintain safety distance is not suitable for practical use. For instance, when the ship turns, the rectangular points around the ship used in the previous study passed very close to obstacles due to optimize $t_{\mathrm{f}}$ to be the shortest.Such a trajectory may shorten time of berthing, but maneuvering by tracking it may cause discomfort and fear to passengers and nearby manned vessels in practical use. Additionally, the trajectory passing near obstacles is not safe to be used as a reference because there is not enough margin to recover from a gust. Hence, to optimize the trajectory planning in a real port, an optimal trajectory should be searched among trajectories that keep a passing distance away from obstacles as a human captain does.}

\Add{To overcome these limitations}, the authors proposed a new collision avoidance algorithm on this study. The \Add{main improvement} of the proposed collision avoidance algorithm \Add{compared to the previous method stated above} as follows:\Add{(i) to handle multiple obstacles and arbitrary shape of obstacle (limitation (a) and (b)), both target berth and other obstacles were represented by polygons.} (ii) to maintain sufficient safety distance to the obstacle, the idea of ship domain was introduced to collision avoidance algorithm \Add{(limitation (c)).}

\Add{The key features of proposed method are: (i) the collision avoidance with sufficient safety distance to obstacle done by minimizing the penalty function $C$ which is the interference of ship domain and obstacle; (ii)} to simplify the algorithm, ship domain was represented by polygon which has vertices at the boundary of domain (hereinafter, these vertices called as ``domain boundary vertices''); \Add{(iii) the ship domain varies it size with ship's speed. This is to enable the ship to maintain distance from obstacle and approach to target berth when the ship reduced its speed; (iv) to generate a favorable gradient of \Cref{eq:objfunc}, penalty function $C$ was expressed as integral of instantaneous interference. A favorable gradient is necessary to lead the update process of distribution of candidate solution to be less interfered. }

\Add{The flow of proposed algorithm avoid are follows. First,  the state vector $\boldsymbol{x}(t)$ for solution $\boldsymbol{X}$ is fed to collision avoidance algorithm. Second, location of domain boundary vertices $P_{j=1\cdots n_{v}}$ are computed by $\boldsymbol{x}(t)$. Third, interference between $P_{j=1\cdots n_{v}}$ and obstacle polygon $S_{\text{obstacle},i=1\cdots n_{\text{obstacle}}}$ are computed by \Cref{eq:smallc} as instantaneous penalty. Finally, sum of $ c_{i,j}(t)$ is derived as penalty function $C$ of solution $\boldsymbol{X}$ and fed to \Cref{eq:objfunc}. Note that proposed collision avoidance algorithm are incorporated as a part of optimization process shown on \Cref{fig:flowchart}. Details of each part of proposed algorithm are shown below.}

\Add{The penalty function $C$ was defined as summation of} the penetration length of domain boundary vertex to the obstacle \Add{for whole duration of berthing trajectory:}
\begin{align}
    \label{eq:penetration_penalty}
        C =\int_{0}^{t_{f}} \sum_{i=n_{\text {obstacle }}} \sum_{j=n_{\text {v}}} c_{i, j}(t) dt \enspace,
\end{align}
where $c_{i,j}(t)$  is instantaneous penalty between $i$- th obstacle polygon $S_{\text {obstacle},i}$ and $j$-th domain boundary vertex $P_{j}$: 
\begin{align}\label{eq:smallc}
        c_{i,j}(t) =
        \begin{cases}
        L_{\text {penet},j}&: P_{j} \in S_{\text {obstacle},i} \\
        0&: P_{j} \notin S_{\text {obstacle},i}
        \end{cases} \enspace.
\end{align}
Here, $n_{\text {v}}$ and $n_{\text {obstacle}}$ represents the total number of domain boundary vertices and obstacles, $L_{\text {penet},i,j}$ is the penetration length of $P_{j}$ to $S_{\text {obstacle},i}$, which computed from the length between $P_{j}$ and nearest edge of $S_{\text{obstacle},i}$ to $P_{j}$.

\paragraph{\Add{Details of ship domain}}
Again, in this study, the ship domain was applied to the collision avoidance algorithm. \Add{Details of ship domain and domain boundary vertices are shown below.} There are two major factors to define the ship domain; the shape of the domain and its size. Regarding the shape of the domain, the shape was designed as an ellipse which has the length ratio of 3:2:1 from midship to the bow, stern, and both sides of a ship, respectively, by referring the research on minimum passing distance in harbor based on inquires to pilots and captains \citep{INOUE1994}. 

In terms of the size of the domain, to approach to the berth, the ship domain has to shrink its size to avoid interference between the ship domain and the berth. Moreover, the ship must slow down when maneuvering in the vicinity of obstacle; it is better to keep a longer passing distance due to the increase of minimum stopping distance with speed. 

Hence, the $P_{i}$'s coordinate $P_{i,x} \ P_{i,y}$ are determined by the semi-major axis $a$, semi-minor axis$b$ and angle $\alpha_{i}$: 
\begin{align}
     (P_{i,x},P_{i,y}) = (a\cos\alpha_{i}, b\sin\alpha_{i})\enspace, 
\end{align}
where
\begin{align}
\begin{aligned}
    a &= L_x(U) + 0.5 L_{pp} \\
    b &= L_y(U) + 0.5 B\enspace.
\end{aligned}
\end{align}
Margin length $L_{x}$ and $L_{y}$ are the function of the resultant velocity $U=\sqrt{u ^{2}+v_{m}^2} \ $. Angle $\alpha_{i}$ is equally spaced, and defined by the ship-fixed coordinate system $o-xy$, with zero in the $x$ direction and positive in the clockwise direction. 
% Moreover, when maneuvering inside a port, including berthing and unberthing, the ship's speed varies greatly compared to navigation in the open sea. Hence, the passing distance to the obstacle should be taken into account the effect of speed. If the speed is low, the stopping distance will be shortened, and it will be safe to approach a fixed obstacle more closely. Accordingly, the ship domain in this study increases its size according to the speed, i.e., the area is larger when the ship is fast, and shrink its area when the speed is low for berthing or entering a narrow fairway. The ship domain that the domain size increased with the quadratic function of speed but tapered off at higher speed was proposed \citet{Wang2016}.
% From the above, the ship domain of this study was defined as a tear-drop shaped ellipse which has the length ratio of 3:2:1 on bow, stern and width direction. The ship domain was simplified as a polygon which has $n_{\text{detect}}$ vertices as domain boundary vertices to reduce computational cost. Despite the original definition of the ship domain is the area which is free from other ship, we extended the definition to the area which is free from obstacle. The coordinates of the proposed domain boundary vertices $P_{i}=(P_{x,i},P_{y,i})$ :
% \begin{align}
%      P_{i} = \boldsymbol{g}\{U, \alpha_{i}\}
% \end{align}
the $x$ direction margin length $L_{x}$ has different lengths $L_{x, L} > L_{x, S}$ in front and behind the hull. The length of the long side of the major axis is switched according to the forward speed's sign:
\begin{align}
    \begin{aligned}
        u \geq 0 :& \begin{cases}
                  L_{x, L} \leftarrow \alpha_{i} \in (-\pi/2, \pi/2) \\
                  L_{x, S} \leftarrow \alpha_{i} \notin (-\pi/2, \pi/2) \\  
        \end{cases}\\
        u < 0 :& \begin{cases}
                  L_{x, L} \leftarrow \alpha_{i} \notin (-\pi/2, \pi/2) \\
                  L_{x, S} \leftarrow \alpha_{i} \in (-\pi/2, \pi/2) \enspace .
                \end{cases}
    \end{aligned}
\end{align}
The $L_{x, L}, \ L_{x, S}, \textrm{and} \ L_{y}$ were constrained to certain maximum and minimum size to prevent the ship domain from becoming too large with increasing speed and to maintain a certain passing distance from obstacles even when the speed is low enough. To make the ship domain of this study applicable to the change of geometry of the port, the max. size of the ship domain was determined by the minimum passage width $W$ in the port. \Add{This is done by following the idea of dynamic ship domain model of \citep{Liu2016} which use track width for domain width on the navigation at restricted channels.} We assumed a distance of 0.25$W$ from the edge of the obstacle in the width direction to pass the head-on situation. Combining with the length ratio of the ellipse's axes, we can obtain: 
\begin{align}\label{eq:max_L}
    \begin{aligned}
        L_{x, \max. L} &=0.75 W-0.5 L_{pp} \\
        L_{x, \max, S} &=0.5 W-0.5 L_{pp} \\
        L_{y,\max} &=0.25 W-0.5 B \enspace.
    \end{aligned}    
\end{align}
On The minimum size, the shape is assumed to be front-to-back symmetric, to maintain sufficient margin on both bow and stern when rotating the ship, which may occur on the terminal phase of berthing and initial phase at confined berth. $L_{y},\min$ was set to ship's width $B$, which is a typical stop point before mooring:
 \begin{align}\label{eq:min_L}
    \begin{aligned}
        L_{x, \min} &=0.25L_{pp} \\
        L_{y,\min}&= B \enspace.
    \end{aligned}    
\end{align}
Finally, The $L_{x, L}, \ L_{x, S}, \textrm{and} \ L_{y}$ express as follows:
\begin{align}
    \begin{aligned}
    &L_{x, L}(U) =&\\
    &\begin{cases}
        L_{x, \max, L} &: U>U_{\max } \\
        L_{x, \min}
        +\frac{\left(L_{x, \max, L}-L_{x,min}\right)\left(U-U_{\min }\right)}{U_{\max }-U_{\min }} &: U_{\min } \leq U \leq U_{\text {max }} \\
        L_{x, \min}  &:  U<U_{\min }
    \end{cases}
    \end{aligned}
\end{align}
\begin{align}
    \begin{aligned}
    &L_{x,S }(U)= &\\
    &\begin{cases}
        L_{{x, \max, S}} &: U>U_{\max } \\
        L_{x, \min}+\frac{\left(L_{{x, \max, S}}-L_{x, \min}\right)\left(U-U_{\min }\right)}{U_{\max }-U_{\min }} &: U_{\min } \leq U \leq U_{\text {max }} \\
        L_{x, \min}  &: U<U_{\min }    
    \end{cases}
\end{aligned}
\end{align}
\begin{align}
    \begin{aligned}
    &L_{y}(U)=& \\
    &\begin{cases}
        L_{y} &: U>U_{\max } \\
        L_{y, \min}+\frac{\left(L_{y, \max}-L_{y, \min }\right)\left(U-U_{\min }\right)}{U_{\max }-U_{\min }} &: U_{\min } \leq U \leq U_{\text {max }} \\
        L_{y, \min} &: U<U_{\min } \enspace,    
    \end{cases}
\end{aligned}
\end{align}
where linearly changing size with speed over the interval of speed $U\in[U_{\min},U_{\max}]$. The limits in this study were set to $[U_{\min},U_{\max}] = [1 \text{kn}, 6 \text{kn}]$.

\begin{figure}
    \centering
    \includegraphics[width=\columnwidth]{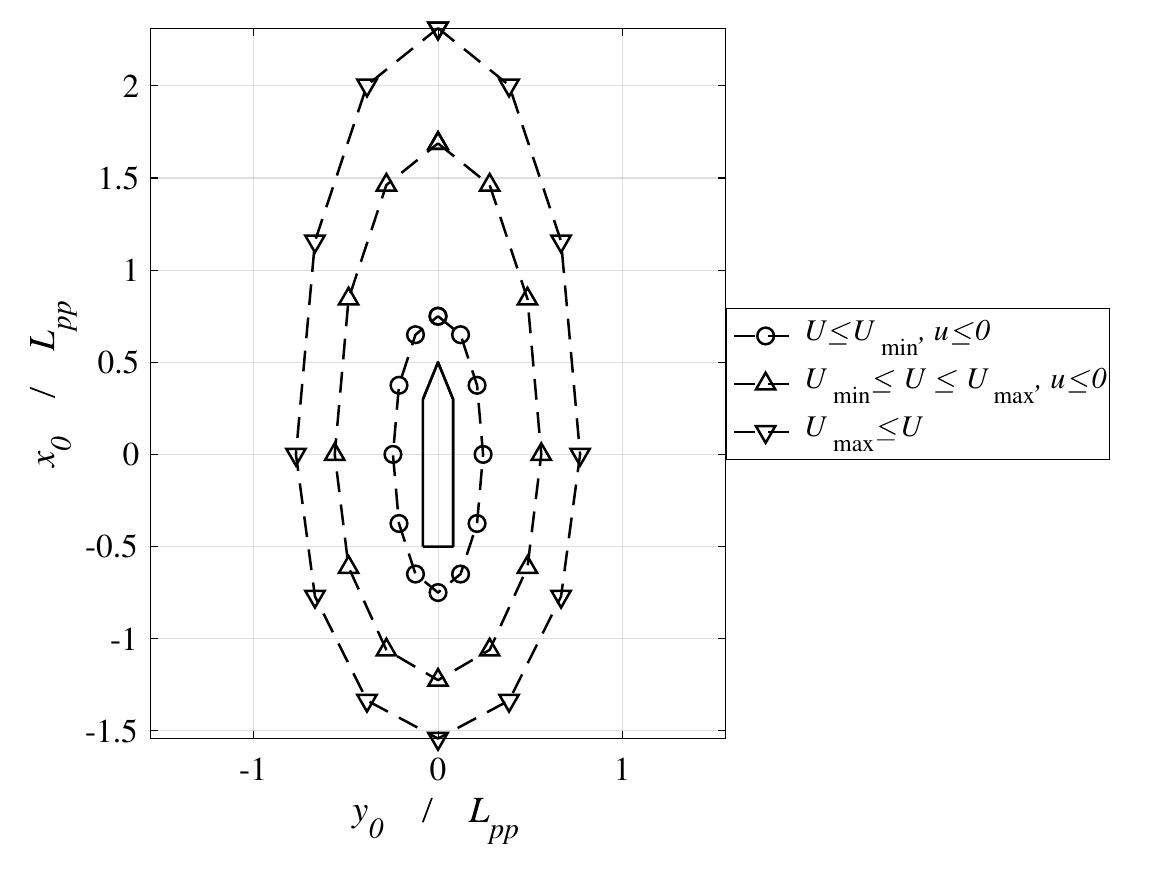}
    \caption{domain boundary vertices on various speed}
    \label{fig:collision_detection_sample}
\end{figure}

\Cref{fig:collision_detection_sample} shows the ship domain and domain boundary vertices for $W=3.08 \ L_{\mathrm{pp}}$, which is the minimum passage width of the port of  \textit{Nanko} shown in \cref{sec:port_geo}. The number of domain boundary vertices $n_{\text{detect}}=13$ in this study.

\subsection{Mathematical Model of Maneuvering}\label{sec:MMG}
In this section, the mathematical model of ship maneuvering is presented. In this study, the Mathematical Maneuvering Group model (MMG model) \citep{Ogawa1978,Yasukawa2015} \EraseTwo{is}\AddTwo{was} applied.\AddTwo{The MMG model is a mathematical model that is considered to be as related as possible to the physical meaning of hydrodynamics \citep{Ogawa1978}, and one of the most common models for ship maneuvering. The MMG model of the subject ship had been validated with both model ship experiments \citep{Maki2020b} and full-scale ship trail \citep{Kobayashi2002}.}

\subsubsection{Subject Ship}
The subject ship of this study is the same as previous \citep{Maki2020}. The subject ship is VLCC M.V. \textit{Esso Osaka}, a single-rudder and single-fixed pitch propeller ship. Although actual \textit{Esso Osaka} does not have any side thruster, to achieve berthing without support from the tug, numerical bow thruster and stern thruster model were added.  \textit{Esso Osaka} was chosen because of its availability of the maneuvering model, while the actual ships operated on the subject ports shown on \cref{sec:port_geo} are passenger ferry and RORO cargo ship. These actual ships operated without the assistance of tug during berthing and unberthing and were equipped with side thrusters and a controllable pitch propeller. Although the size was a model scale in the previous study \citep{Maki2020}, it was expanded to actual ship size. To be equivalent to the ships actually operated on the subject ports, the ship size was shrunken to $L_{\textrm{pp}}=150$ m while Original \textit{Esso Osaka} is $L_{\textrm{pp}}=325$ m.
The principal particulars of the subject ship are shown in \Cref{tab:pp}.

\begin{table}
    \centering
    \caption{Principal particulars of subject ship}
    \begin{tabular}{ll}
    \hline Item & Value \\
    \hline Length between perpendicular: $L_{\mathrm{pp}}$ (m) & $150.0$ \\
    Ship breadth: $B$ (m) & $24.46$ \\
    Ship draft: $d$ (m)& $10.06$ \\
    Diameter of propeller: $D_{\mathrm{p}}$ (m)& 4.20 \\
    Area of Rudder: $A_{\mathrm{R}}$ ($\mathrm{m}^{2})$ & 26.58\\
    Diameter of bow thruster: $D_{\mathrm{BT}}$: (m)& 2.5 \\
    Diameter of stern thruster: $D_{\mathrm{ST}}$: (m)& 2.5 \\
    Mass: $m$ (ton) & 31,412\\
    Longitudinal center of gravity: $x_{\mathrm{G}}$ (m)& 4.75\\
    Transverse projected area:$A_{\mathrm{T}}$($\mathrm{m}^{2})$& 213.55\\
    Lateral projected area: $A_{\mathrm{L}}$($\mathrm{m}^{2})$& 1150.50\\
    Block coefficient: $C_{\mathrm{b}}$ & 0.831 \\
    \hline
    \end{tabular}
    
    \label{tab:pp}
\end{table}

\subsubsection{MMG model}
The 3-DoF equation of motion of MMG model is express as follows:
\begin{equation}
    \begin{split}
    (m + m_x) \dot{u} - (m + m_y) v_m r -x_{G}mr^{2} &= F_{x}\\
    (m + m_x) \dot{v}_m - (m + m_y) u r +x_{G}m\dot{r} &= F_{y}\\
    (I_{zz} + J_{zz}+x_{G}^{2}m) \dot{r} - x_G m (\dot{v}_m + ur) &= M_{N}
    \end{split}
\end{equation}
with
\begin{equation}
    \begin{split}
    F_{x} &= X_H + X_P + X_R + X_A \\
    F_{y} &= Y_H + Y_P + Y_R +  Y_{SS} + Y_A\\
    M_{N} &= N_H + N_P + N_R + N_{SS} + N_A\enspace.
 \label{eq:MMGdynamics}
    \end{split}
\end{equation}
Since the MMG model decomposes the hydrodynamic force acting on the ship to sub-model for major components consisting of the ship, the right-hand side of \cref{eq:MMGdynamics} represents the force or moment acting on each component of the ship; the subscripts H, P, R, A and SS denote the hull, the propeller, the rudder, the external forces by wind and the side thrusters, respectively. Each sub-model contains model coefficients to estimate force form input $\boldsymbol{u}$ and state $\boldsymbol{x}$. The model coefficients, also known as hydrodynamic coefficients, are vessel-specific, and they may vary by the scale of the ship; due to the large discrepancy of the Reynolds number between a scale model for towing tank experiments and the actual ship.
However, in this study, model coefficients obtained by scale model test (\citep{Hachii2004, Kobayashi2002} were used, because the values for $L_{\text{pp}}=150$ m of \textit{Esso Osaka} are unknown.

Choice of sub-model follows the previous research \citep{Maki2020}, with several new additions. For the force acting on the hull, Yoshimura's unified model \citep{Yoshimura2009} was applied. For the force induced by the propeller, the model is switched based on the operating condition of the propeller. For the propeller forward rotation ($n_{\mathrm{p}}\geq 0$), thrust force $X_{p}$ expressed by standard MMG model \citep{Yasukawa2015}.
For the propeller reversal condition ($n_{\mathrm{p}}< 0$), propeller force and moment are modeled by a polynomial equation based on towing tank experiment \citep{Hachii2004, Ueno2001}. For the force induced by rudder, standrd MMG model was applied for forwarding ($u\geq 0)$ condition and 4th quadrant operation($u<0, n_{\mathbf{p}}<0$), Kitagawa's model \citep{Kitagawa2015} for 3rd quadrant operation ($u<0, n_{\mathbf{p}}\geq0$). Refer \citet{Miyauchi20201} for the detail of the MMG model used in this paper.

Regarding the external force induced by wind disturbance, Fujiwara's regression formulae \citep{Fujiwara1998} was used to estimate the wind pressure coefficients:
\begin{align}
    \begin{aligned}
        X_{A} &= (1/2)\rho_{A}U_{A}^{2}A_{T}\cdot C_{X} \\
        Y_{A} &= (1/2)\rho_{A}U_{A}^{2}A_{L}\cdot C_{Y} \\
        N_{A} &= (1/2)\rho_{A}U_{A}^{2}A_{L}L_{OA}\cdot C_{N} \enspace,
    \end{aligned}
\end{align}
where
\begin{align}
    \begin{aligned}
        C_{X} =& X_{0}+X_{1} \cos (2\pi - \gamma_{A})+X_{3} \cos 3 (2\pi - \gamma_{A}) \\
              &+X_{5} \cos 5 (2\pi - \gamma_{A}) \\
        C_{Y} =& Y_{1} \sin (2\pi - \gamma_{A})+Y_{3} \sin 3 (2\pi - \gamma_{A}) \\
              &+Y_{5} \sin 5 (2\pi - \gamma_{A}) \\
        C_{N} =& N_{1} \sin (2\pi - \gamma_{A})+N_{2} \sin 2 (2\pi - \gamma_{A}) \\
              &+N_{3} \sin 3 (2\pi - \gamma_{A}) \enspace,
    \end{aligned}
\end{align}
$\rho_{A}$ is the density of air, $A_{T},\ A_{L}, L_{OA}$ are the transverse projected area, the lateral projected area, and the overall length of the ship, respectively. $X_{i}, \ Y_{i}, N_{i}$ are coefficients to express wind pressure coefficients derived by the regression formulae, which use geometric parameters of the ship as explanatory variables and based on wind tunnel experiment data of numerous scaled ship models.

For the side thruster model, Kobayashi's model \citep{Kobayashi1988, Kobayashi1990} was used, which follows
\begin{align}
    \begin{aligned}
        X_{ST} =&0 \\
        Y_{SS} =& (1+a_{YSB}\cdot |\rm{Fr}|)\cdot \it T_{BT} \\
        &+ (1+a_{YST}\cdot |\rm{Fr}|)\cdot \it T_{ST}  \\
        N_{SS} =& (1+a_{NSB}\cdot |\rm{Fr}|)\cdot \it T_{BT} \cdot x_{BT} \\ 
        & +  (1+a_{NST}\cdot |\rm{Fr}|)\cdot \it T_{ST} \cdot x_{ST} \\
        T_{BT} =& \rho D_{BT}^4 n_{BT}^{2}K_{TBT} \\
        T_{ST} =& \rho D_{ST}^4 n_{ST}^{2}K_{TST}\enspace .
\end{aligned}
\end{align}
Here, $a_{YBT}, a_{NBT}, a_{YST}, a_{NST}$ are the coefficients expressing the interaction between hull and thruster. $x_{BT}, x_{ST}$ are the longitudinal location of side thruster. In addition to the interaction effect, the thrust of side thrusters are zero when the longitudinal speed is larger than the threshold: $|u| > u_{\text{threshold}}$. The value of threshold $u_{\text{threshold}} = 2.5722$ (m/s)  is equivalent to 5 knots in this study.

\section{Computation Results}
\subsection{Port Geometry}\label{sec:port_geo}
\begin{figure}
 \begin{minipage}[c]{\linewidth}
  \centering
  \includegraphics[keepaspectratio, width=0.8\columnwidth]{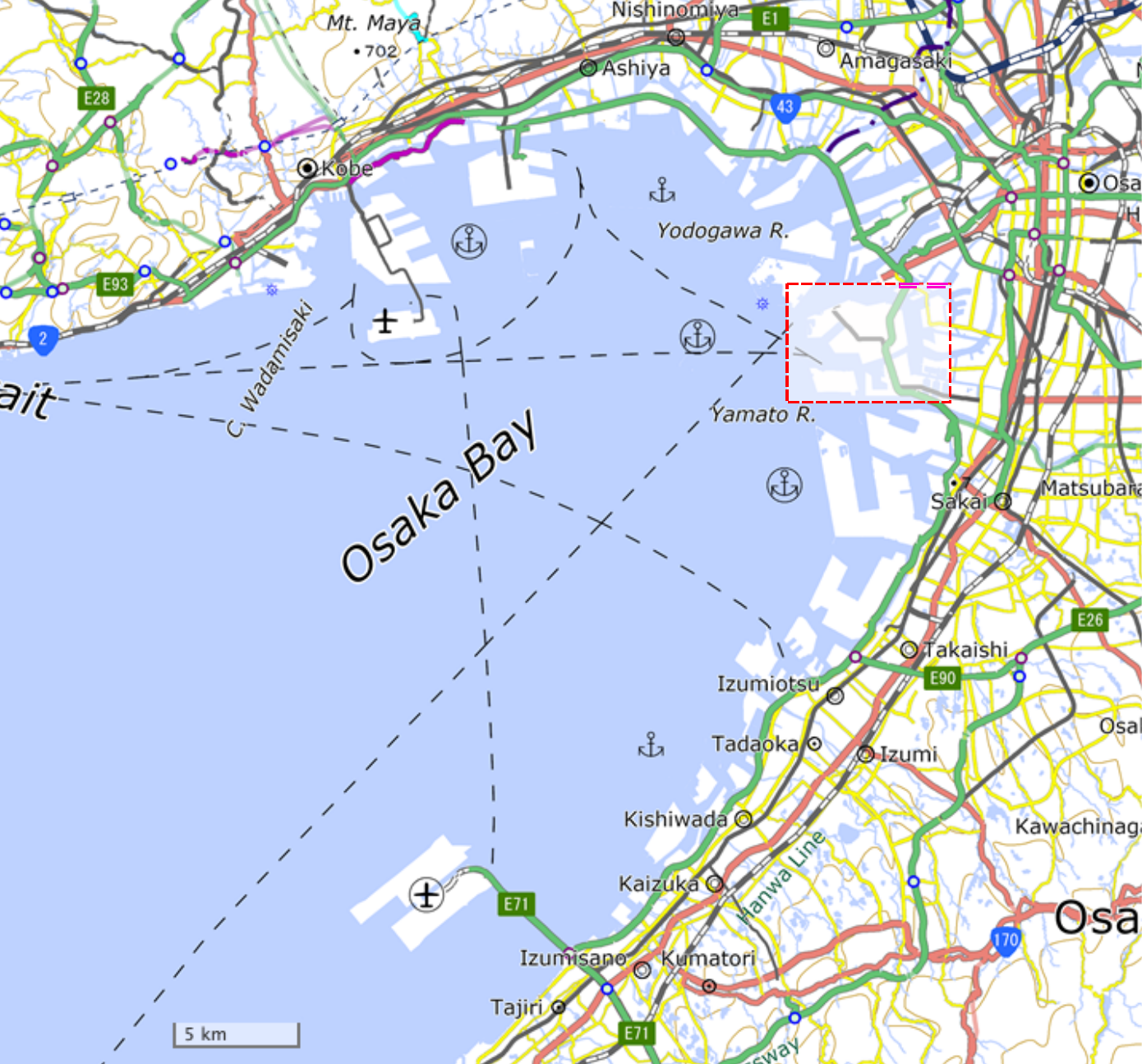}\\
  \subcaption{Far view; red dashed line rectangular shows the area of the figure below}\label{fig:nanko_farview}
 \end{minipage}
 \begin{minipage}[c]{\linewidth}
  \centering
  \includegraphics[keepaspectratio,width=\columnwidth]{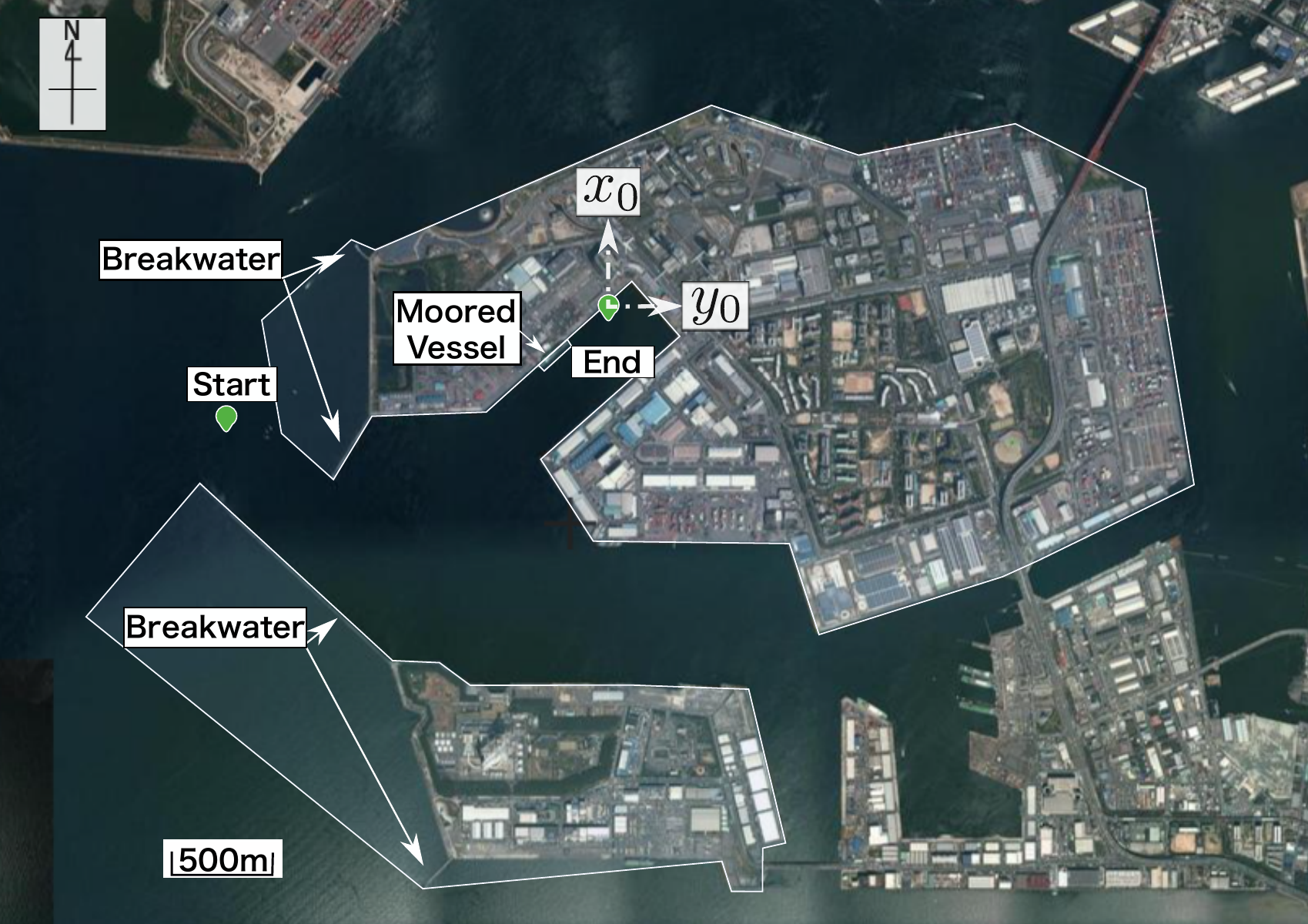}\\
  \subcaption{Near view. polygons with white line represent the obstacles which modeled in the computation.}\label{fig:nanko_sat}
 \end{minipage}
 \caption{Port of  \textit{Nanko}. Created by Authors based on 1:1,000,000 INTERNATIONAL MAP data and seamless photo data (Geospatial Information Authority of Japan) (http://maps.gsi.go.jp/)}\label{fig:Nanko_photo}
\end{figure}

\begin{figure}
 \begin{minipage}[c]{\linewidth}
  \centering
  \includegraphics[keepaspectratio, width=0.8\columnwidth]{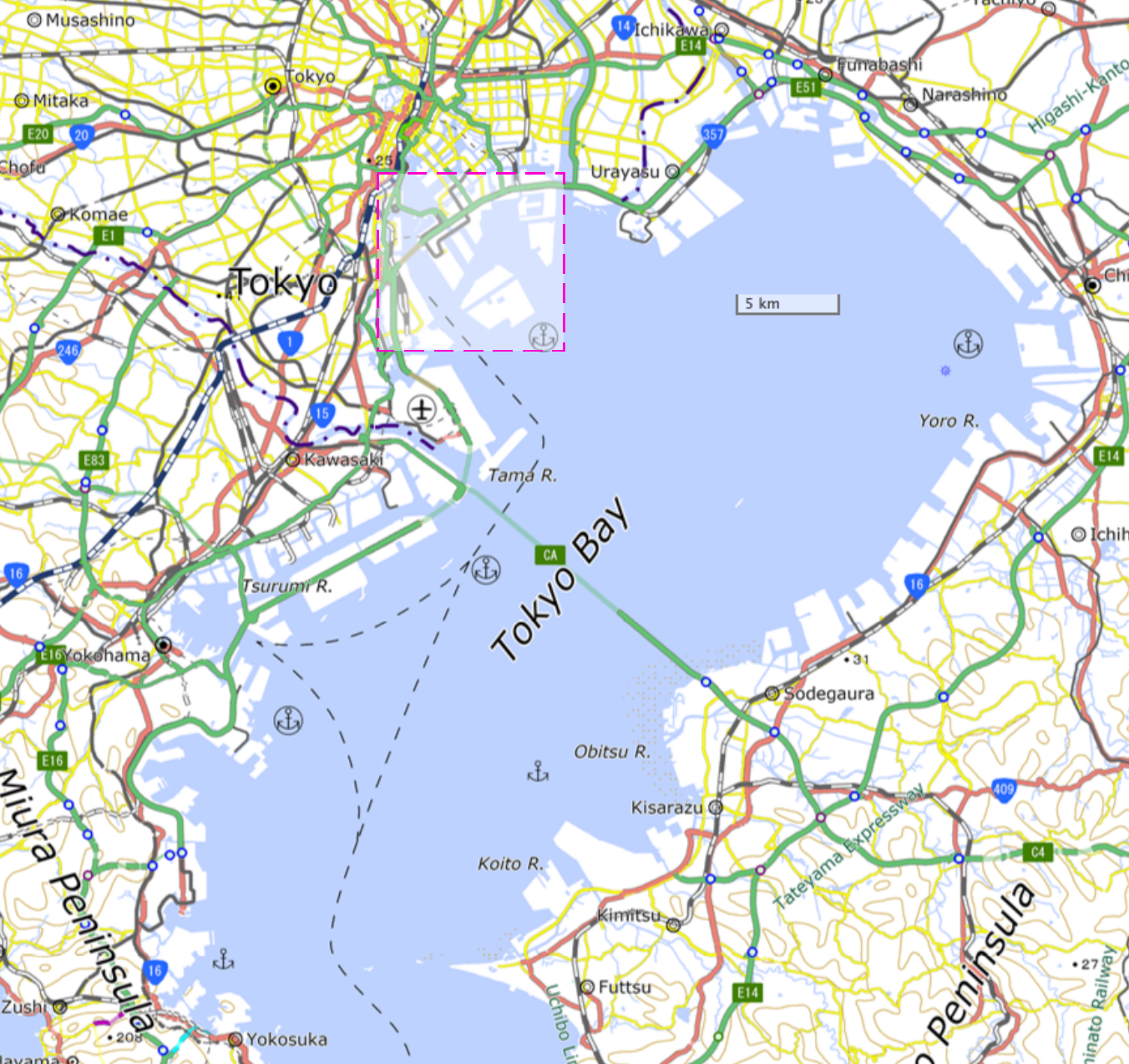}\\
  \subcaption{Far view; red dashed line rectangular shows the area of the figure below}\label{fig:ariake_farview}
 \end{minipage}
 \begin{minipage}[c]{\linewidth}
  \centering
  \includegraphics[keepaspectratio,width=0.83\columnwidth]{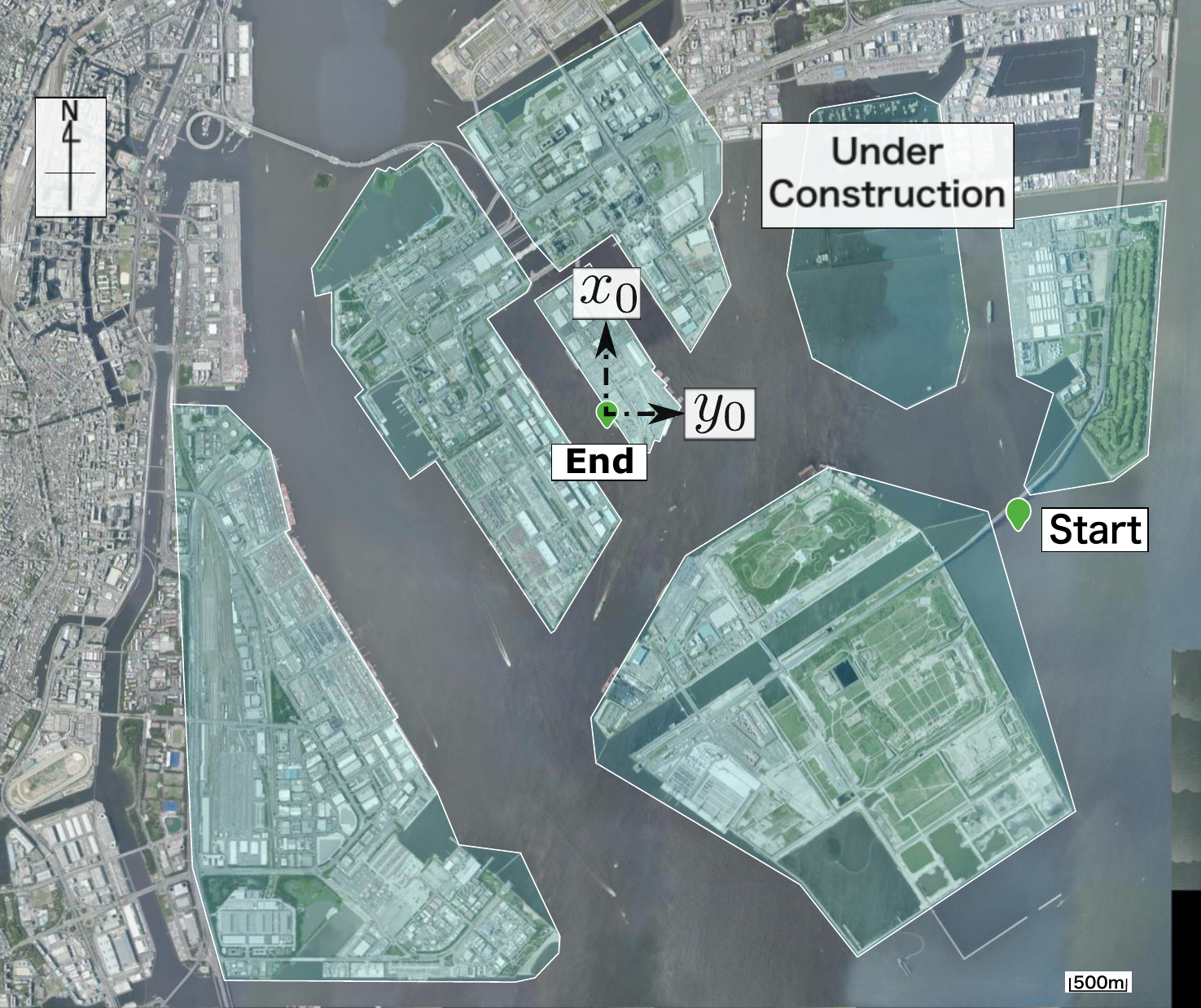}\\
  \subcaption{Near view. polygons with white line represent the obstacles which modeled in the computation.}\label{fig:ariake_sat}
 \end{minipage}
 \caption{Port of  \textit{Ariake}. Created by Authors based on 1:1,000,000 INTERNATIONAL MAP data and seamless photo data (Geospatial Information Authority of Japan) (http://maps.gsi.go.jp/)}\label{fig:Ariake_photo}
\end{figure}

Two ports in Japanese waters were selected for computation; Port of \textit{Nanko} at Osaka bay and the port of \textit{Ariake} at Tokyo bay.  \textit{Nanko} is the primary port in this study, while  \textit{Ariake} is alternate to verify the applicability to arbitrary port geometry of the present method. Overview of ports is shown in \Cref{fig:Nanko_photo} and \Cref{fig:Ariake_photo}. These figures also show the start points and the endpoints of berthing. \Cref{tab:nanko_condition} and \Cref{tab:ariake_condition} show the detailed value of $\boldsymbol{x}_{\mathrm{int}}$ and $\boldsymbol{x}_{\mathrm{des}}$ of computation. As stated in the \cref{sec:method}, the origin of $O_{0}-x_{0}y_{0}$ coordinates system is set to the endpoint for berthing and the start point for unberthing. The start point of berthing was set just outside the breakwaters of the port. This is because, outside the breakwater, the spatial constrain is not severe as inside; hence the appropriate collision avoidance method might differ. The endpoint was set at the point which perpendicular distance is 40 m from berth wall; Approximately 1.0 $B$ from the broadside of the ship. The berthing condition is the starboard side, head out at  \textit{Nanko} and port side, head out at  \textit{Ariake}. Note that the obstacle polygons include certain water surfaces at the vicinity of the breakwater as a restricted area.  

Both ports have characteristics to make berthing and unberthing difficult: dead end, confined berth which ship necessary to turn 180 degrees, require several course change during navigation, and obstacle landfills that exclude the intended berth. The case of  \textit{Nanko} also has a moored vessel near the berthing point as an obstacle. Generally, the ship must comply with the international and local traffic rules; the ship must navigate in the designated passage in a port (for the regulation on the port of  \textit{Nanko}, see \citet{nanko2020}). However, in this study, regulatory requirements such as passage and fairway are not considered.

The minimum passage width $W$, which is the design parameter of domain boundary vertices, are 3.08 $L_{\mathrm{pp}}$ for  \textit{Nanko} and 2.40 $L_{\mathrm{pp}}$ for  \textit{Ariake}, respectively.

\begin{table*}
    \centering
    \caption{Computation condition of port of  \textit{Nanko}}
    \begin{tabular}{cc}
        \hline Parameter & Value \\
        \hline \hline Condition & Berthing \\
        $\boldsymbol{x}_{\mathrm{int}}$ & $(-598.7 \ \mathrm{m}, 8.0 \ \mathrm{kn}, -1845.2 \ \mathrm{m}, 0.0 \ \mathrm{m/s}, 132 ^{\circ}, 0.0^{\circ}/\mathrm{s})^\mathsf{T} $ \\
        $\boldsymbol{x}_{\mathrm{des}}$ &  $(0.0 \mathrm{m}, 0.0 \mathrm{kn}, 0.0 \mathrm{m}, 0.0 \mathrm{m/s}, 227 ^{\circ}, 0.0^{\circ}/\mathrm{s})^\mathsf{T}  $ \\
        \hline
        Condition & UnBerthing \\
        $\boldsymbol{x}_{\mathrm{int}}$ & $(0.0\ \mathrm{m}, 0.0 \ \mathrm{kn}, 0.0 \ \mathrm{m}, 0.0 \ \mathrm{m/s}, 227 ^{\circ}, 0.0^{\circ}/\mathrm{s})^\mathsf{T}  $ \\
        $\boldsymbol{x}_{\mathrm{des}}$ & $(-598.7 \ \mathrm{m}, 6.0 \ \mathrm{kn}, -1845.2 \ \mathrm{m}, 0.0 \ \mathrm{m/s}, 312 ^{\circ}, 0.0^{\circ}/\mathrm{s})^\mathsf{T} $\\
        \hline
    \end{tabular}
    \label{tab:nanko_condition}
\end{table*}

\begin{table*}
    \centering
    \caption{Computation condition of port of  \textit{Ariake}}
    \begin{tabular}{cc}
        \hline Parameter & Value \\
        \hline \hline Condition & Berthing \\
        $\boldsymbol{x}_{\mathrm{int}}$ & $(-891.3 \ \mathrm{m}, 8.0 \ \mathrm{kn}, 3317.5 \ \mathrm{m}, 0.0 \ \mathrm{m/s}, 326 ^{\circ}, 0.0^{\circ}/\mathrm{s})^\mathsf{T} $ \\
        $\boldsymbol{x}_{\mathrm{des}}$ &  $(0.0 \mathrm{m}, 0.0 \mathrm{kn}, 0.0 \mathrm{m}, 0.0 \mathrm{m/s}, 146 ^{\circ}, 0.0^{\circ}/\mathrm{s})^\mathsf{T}  $ \\
        \hline
        Condition & UnBerthing \\
        $\boldsymbol{x}_{\mathrm{int}}$ & $(0.0\ \mathrm{m}, 0.0 \ \mathrm{kn}, 0.0 \ \mathrm{m}, 0.0 \ \mathrm{m/s}, 146 ^{\circ}, 0.0^{\circ}/\mathrm{s})^\mathsf{T}  $ \\
        $\boldsymbol{x}_{\mathrm{des}}$ & $(-891.3 \ \mathrm{m}, 6.0 \ \mathrm{kn}, -3317.5 \ \mathrm{m}, 0.0 \ \mathrm{m/s}, 146 ^{\circ}, 0.0^{\circ}/\mathrm{s})^\mathsf{T} $\\
        \hline
    \end{tabular}
    \label{tab:ariake_condition}
\end{table*}

\subsection{Berthing and Unberthing Results with proposed ship domain and collision avoidance algorithm}\label{sec:nnk_nowind}
\begin{figure}
    \centering
    \includegraphics[width=\columnwidth]{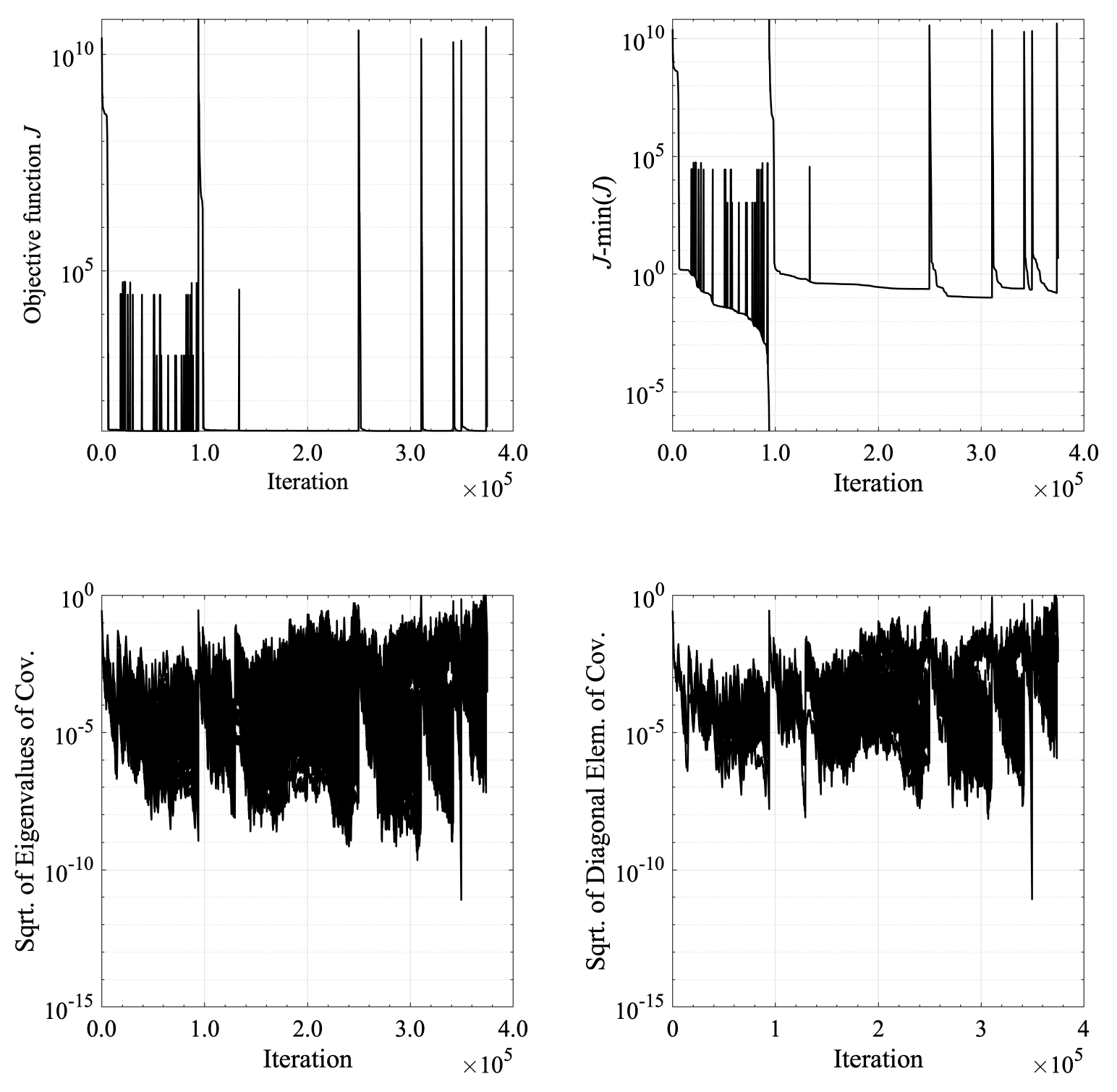}
    \caption{Optimization process of CMA-ES.}
    \label{fig:convergence}
\end{figure}

The trajectory planning of berthing and unberthing are conducted with newly introduced objective function including collision avoidance to the polygonal constraints on the port of  \textit{Nanko}. The convergence process in the optimization is shown in \Cref{fig:convergence}. The upper-left shows the best objective function $J$ at each iteration. The upper-right shows the difference of $J$ at each iteration and the minimum value of $J$ through the optimization process. The lower-left and lower-right show the square root of the eigenvalue of the covariance matrix and the square root of the diagonal elements of the covariance matrix during the process, respectively. As stated in \citet{Maki2020b}, the square root of each eigenvalue represents the axis length of the equidensity ellipsoid, whereas the square root of each diagonal element stands for the standard deviation in each coordinate. From the upper two sub-figures on \Cref{fig:convergence}, $J$ shows impulse-like increases caused by the restart of CMA-ES. By using the restart strategy, CMA-ES let the $J$ converge to several different local minima and choose the best solution from those. In the case of \Cref{fig:convergence}, the optimum solution is obtained at the 94000th iteration. In general, the convergence speed will vary by the computation conditions. However, the maximum iteration number of optimization processes was set to $3\times10^5$ based on the convergence of \Cref{fig:convergence}, which occur several restarting and obtain local minima. The computation was conducted on the workstation equipped with Intel Xeon Gold 6248R for CPU as a serial computation. The computation time took few days to reach maximum iteration number $3\times10^5$.

Since CMA-ES is based on a stochastic method, results will differ by the random seed. \Cref{tab:statistics} show the $t_{\mathrm{f}}$, $J$ and $C$ of 10 independent trials on berthing computation of  \textit{Nanko} without wind disturbance. From the table, we can find the values of $C$ are zero for all 10 trails. This means, not only collision occurred, but also berthing was conducted with sufficient distance to obstacle by maintaining ship domain. With the newly introduced collision avoidance algorithm, collision avoidance to the obstacle of realistic port geometry was successfully achieved.

\begin{table*}
    \centering
        \caption{Statistical analysis data for $J$ and $t_{\mathrm{f}}$in 10 trials}: Minimum, Maximum, Quartiles, Mean, and Standard deviation.
        \begin{tabular}{cccccccc}
            \hline & $\min$ & $\max$ & $Q_{1}$ & $Q_{2}$ & $Q_{3}$ & $\mathrm{Mean}$ & $\sigma$ \\
            \hline$t_{\mathrm{f}}$ & $1081.6$ & $1177.1$ &$
   1086.5$ & $1098.2$&$1103.4$& $1106.0$ & $29.3$  \\
            $J$ & $18.077$ & $19.673$ &$18.159$&$18.355
$&$ 18.407$& $18.484$ & $0.490$  \\
            $ w_{\mathrm{C}}\cdot C$ & $0.0$ & $0.0$ & $0.0$ & $0.0$ & $0.0$ & $0.0$ & $0.0$ \\
            \hline
        \end{tabular}
    \label{tab:statistics}
\end{table*}

\begin{figure*}
    \centering
    \includegraphics[width=\linewidth]{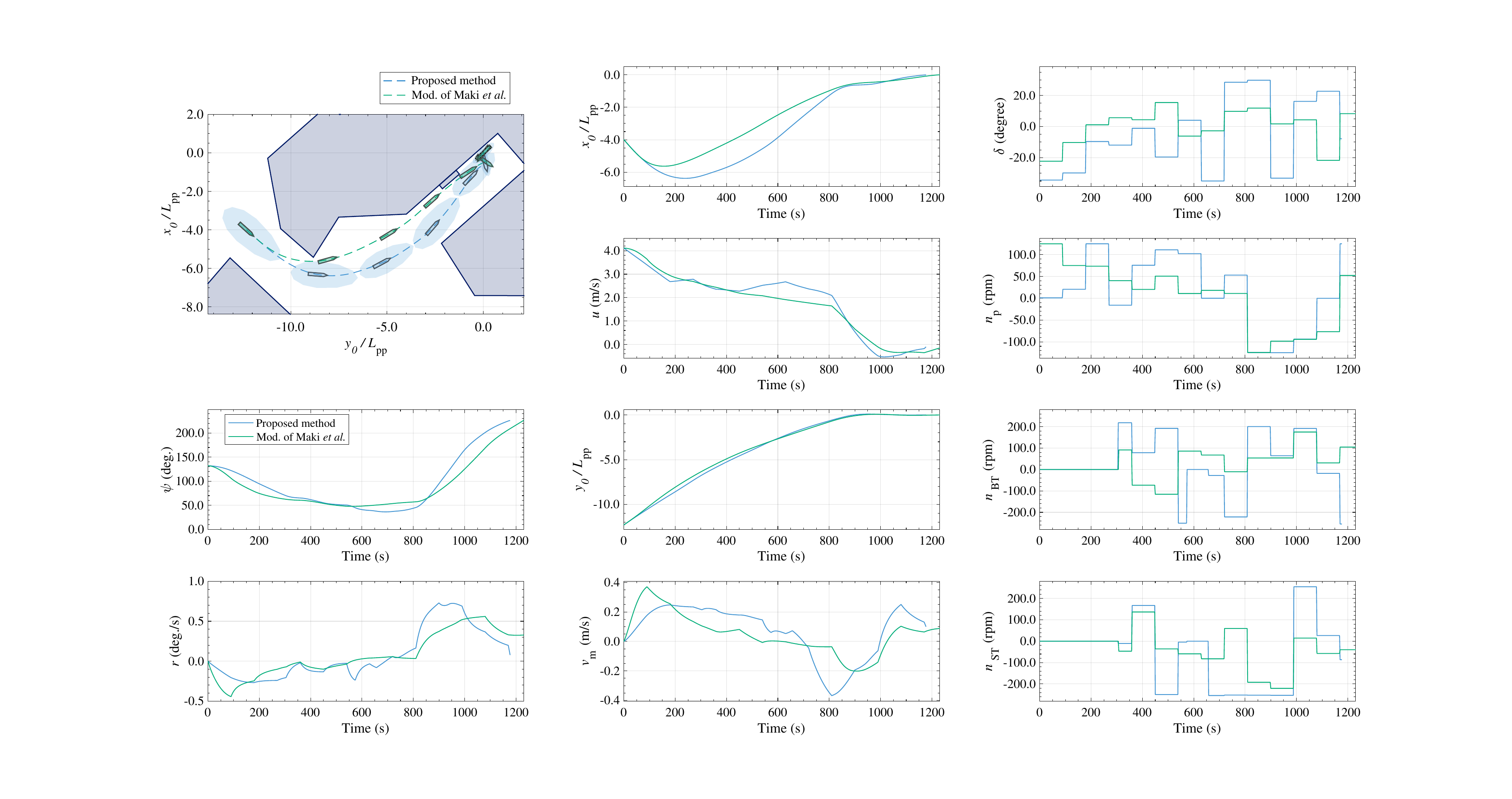}
    \caption{Optimal control inputs and achieved state  of berthing on port of  \textit{Nanko}. Ellipses shown in light blue on figure of trajectory at upper left are the ship domains used in collision avoidance algorithm. \AddTwo{Comparison computation by the previous study \citep{Maki2020} with modification is also shown. The proposed method maintain distance from obstacle when turning while previous method passed without safety distance.}}
    \label{fig:nnk_in_whole}
\end{figure*}

\Cref{fig:nnk_in_whole} shows the obtained control input $\boldsymbol{u}(t)$ and state  $\boldsymbol{x}(t)$ achieved by the optimal $\boldsymbol{u}(t)$ on berthing of  \textit{Nanko} without wind disturbance; $U_{T}=0$. Note that if $u > u_{\text{threshold}}$, which means thrust of side thrusters are zero, shown as $n_{BT}=0$ and $n_{ST}=0$ in the figure. 

\AddTwo{Additionally, a comparison computation with the previous method \citep{Maki2020} is shown in the figure to show the effectiveness of the proposed method to maintain safety. In the previous study, circumscribed rectangular without safety distance was used to represent the ship's hull. For the comparison computation, the objective function of \citet{Maki2020} was modified by introducing the collection vector $\boldsymbol{w}_{\textrm{dim}}$ and $w_{c}=100$ was used.}

\begin{figure*}
    \centering
    \includegraphics[width=\linewidth]{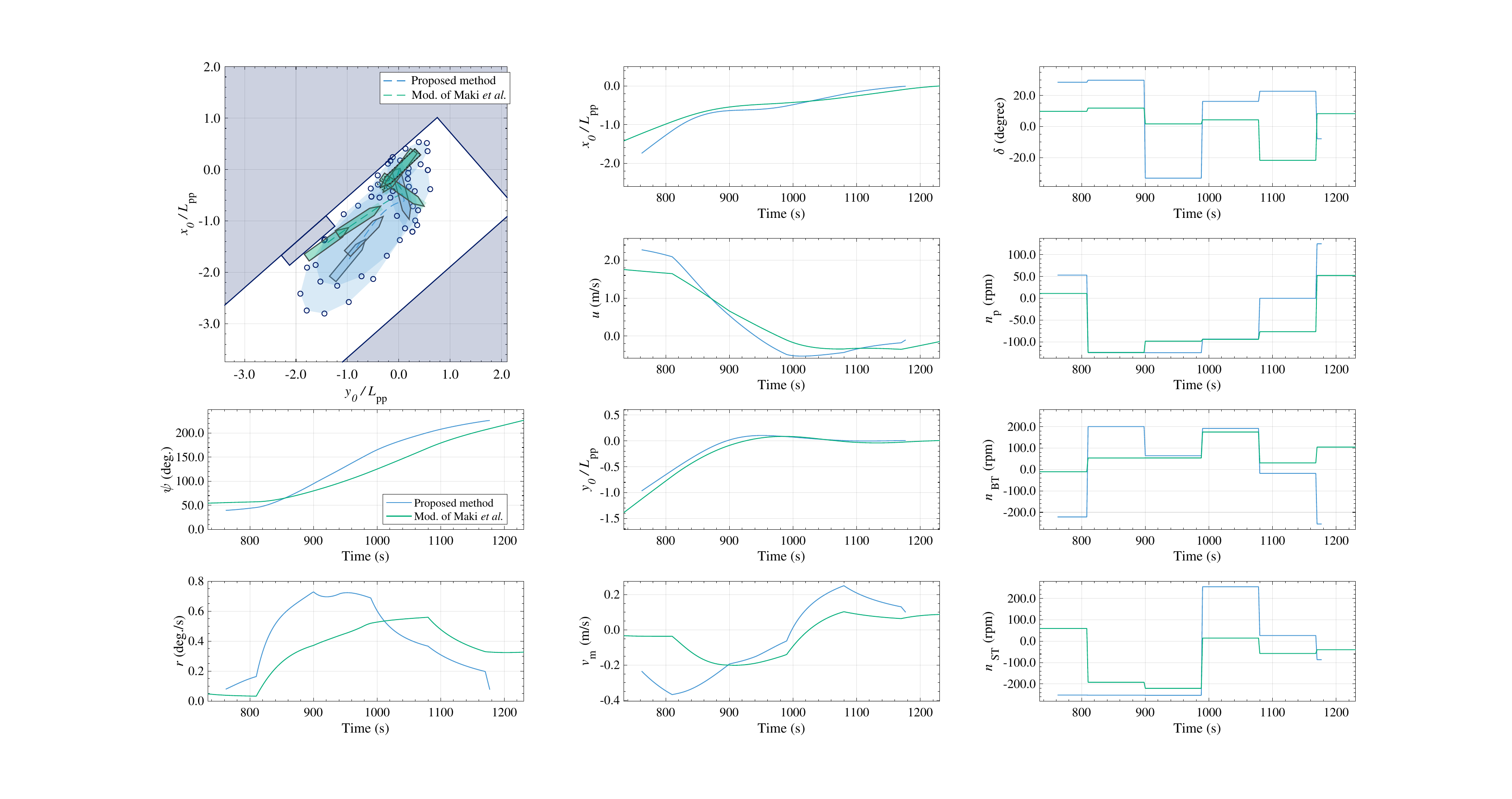}
    \caption{Terminal phase of Fig.\ref{fig:nnk_in_whole} which of the Euclidean distance to berthing point is smaller than 2 $L_{\mathrm{pp}}$. Blue circles on the upper left figure are the collision avoidance points.\AddTwo{ The proposed method passed the moored vessel with distance while previous method passed without safety margin.}}
    \label{fig:nnk_in_near}
\end{figure*}

\EraseTwo{The figure} \Cref{fig:nnk_in_whole} indicates that the optimization on trajectory planning with the new collision avoidance algorithm works well; the ship has reached the berthing point while avoiding the static obstacles with sufficient distance. \Cref{fig:nnk_in_near} shows the terminal phase of the same berthing result with \Cref{fig:nnk_in_whole}; Trajectory and control input are shown only for times when the Euclidean distance between midship and berthing point is less than 2 $L_{\mathrm{pp}}$. In \Cref{fig:nnk_in_whole} and \ref{fig:nnk_in_near}, the light blue ellipses around the ship are ship domain, which shown from $t=0  to \ t=t_{\mathrm{f}}$ at the intervals of 200 seconds. \Cref{fig:nnk_in_near} also shows the domain boundary vertices themselves. As shown in \citet{Maki2020b}, the subject ship tended to turn starboard while propeller reversal. It seems the CMA-ES found the solution which uses the characteristic on the reverse maneuvering of the subject ship to stop within the minimum time.

\AddTwo{\Cref{fig:nnk_in_whole,fig:nnk_in_near} also show the result of comparison computation with the previous method; however,  the previous method shows an inappropriate path. The ship got too close to the obstacle when turning the vicinity of the obstacle (\Cref{fig:nnk_in_whole}) and passing by the moored vessel (\Cref{fig:nnk_in_near}). It should be emphasized that the previous method obtained a collision-free trajectory, which means obtained trajectory is optimal from an optimization scheme standpoint; however, it is insufficient for practical use. In contrast, the proposed method can obtain a trajectory with a certain safety distance suitable for practical use.}


\begin{figure*}
    \centering
    \includegraphics[width=\linewidth]{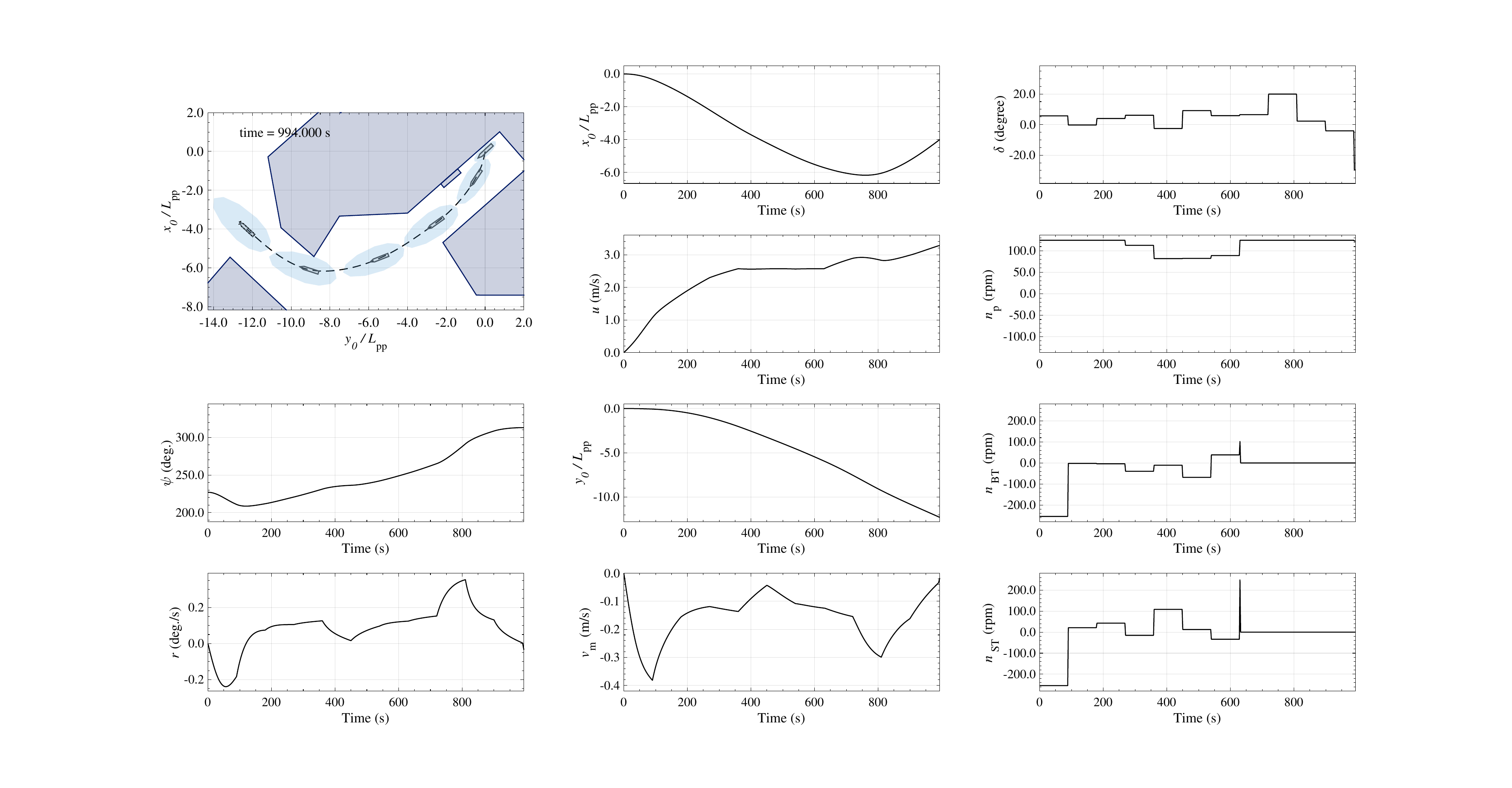}
    \caption{Optimal control inputs and achieved state  of unberthing on port of  \textit{Nanko}. Ellipse shown in light blue on figure of trajectory at upper left is the ship domain used in collision avoidance algorithm. }
    \label{fig:nnk_out_whole}
\end{figure*}

Fig.\ref{fig:nnk_out_whole} shows the result of the unberthing computation. Same as berthing, CMA-ES with the proposed collision avoidance algorithm can obtain resalable trajectory and control inputs. 

Table\ref{tab:results} shows the difference between the desired and the obtained end state: $ \boldsymbol{x}_{\mathrm{des}}-\boldsymbol{x}(t_{\mathrm{f}})$ and $t_{\mathrm{f}}$ for both of berthing and unberthing. For the berthing, $\boldsymbol{x}(t_{\mathrm{f}})$ just fits within the $\boldsymbol{x}_{\text{tol}}$. This is because objective function did not change when  $\boldsymbol{x}(t_{\mathrm{f}})$ was within the range of $\boldsymbol{x}_{\text{tol}}$, while CMA-ES tried to optimize by shortening the $t_{\mathrm{f}}$. Meanwhile, in the unberthing computation, $x_{\mathrm{des},2}-x_{2}(t_{\mathrm{f}}) < 0$ which means the $u(t_{\mathrm{f}})$ is faster than the desired exit speed, while other state satisfy $\boldsymbol{x}_{\text{tol}}$.

\begin{table*}
    \centering
    \caption{Computation results of berthing and unberthing.}
    \begin{tabular}{cc}
        \hline
        parameter & value \\
        \hline\hline
        Condition &  \textit{Nanko}, berthing, $U_{T}=0$\\
        $\boldsymbol{x}_{\mathrm{des}, i}-\boldsymbol{x}_{i}(t_{\mathrm{f}})$ & $(1.0 \mathrm{m}, 0.1  \mathrm{m/s}, -1.0 \mathrm{m}, -0.1 \mathrm{m/s}, 1.0 ^{\circ}, -0.0764^{\circ}/\mathrm{s})^\mathsf{T}  $  \\
        $t_{\mathrm{f}} (s)$ & 1177.1 \\
        \hline
        Condition &  \textit{Nanko}, unberthing, $U_{T}=0$\\
        $\boldsymbol{x}_{\mathrm{des}, i}-\boldsymbol{x}_{i}(t_{\mathrm{f}})$ & $(1.0 \mathrm{m}, -0.2  \mathrm{m/s}, -1.0 \mathrm{m}, 0.02 \mathrm{m/s}, -1.0 ^{\circ}, -0.0343^{\circ}/\mathrm{s})^\mathsf{T}  $  \\
        $t_{\mathrm{f}} (s)$ & 994.0 \\
        \hline
    \end{tabular}
    \label{tab:results}
\end{table*}

\subsection{Result with Wind Disturbance}\label{sec:nnk_wind}
Generally, during the berthing and unberthing, the wind disturbance makes it difficult to control a vessel. Hence, the effect of wind must be included in the optimization of trajectory planning. Moreover, the wind has a critical effect on control a ship at the terminal phase of berthing because the thrusts of actuators are kept low, and the distance to berth is small while the ship drifted by the wind. Obviously, the wind speed and direction are unsteady, and instantaneous velocity can only be described by stochastic approach. Hence, it is not appropriate to implement unsteady wind directly to optimization because CMA-ES would optimize while knowing the time series of instantaneous wind due to its iterative process. Nevertheless, including \textit{steady} wind disturbance is meaningful because considering severe conditions by wind disturbance, the obtained trajectory will be a more generous one that takes into account the limits of actuators. In addition to that, by increasing the wind speed, we can find the nominal wind speed limit of a certain ship and spatial constraints of the port. 

\begin{figure*}
    \centering
    \includegraphics[width=\linewidth]{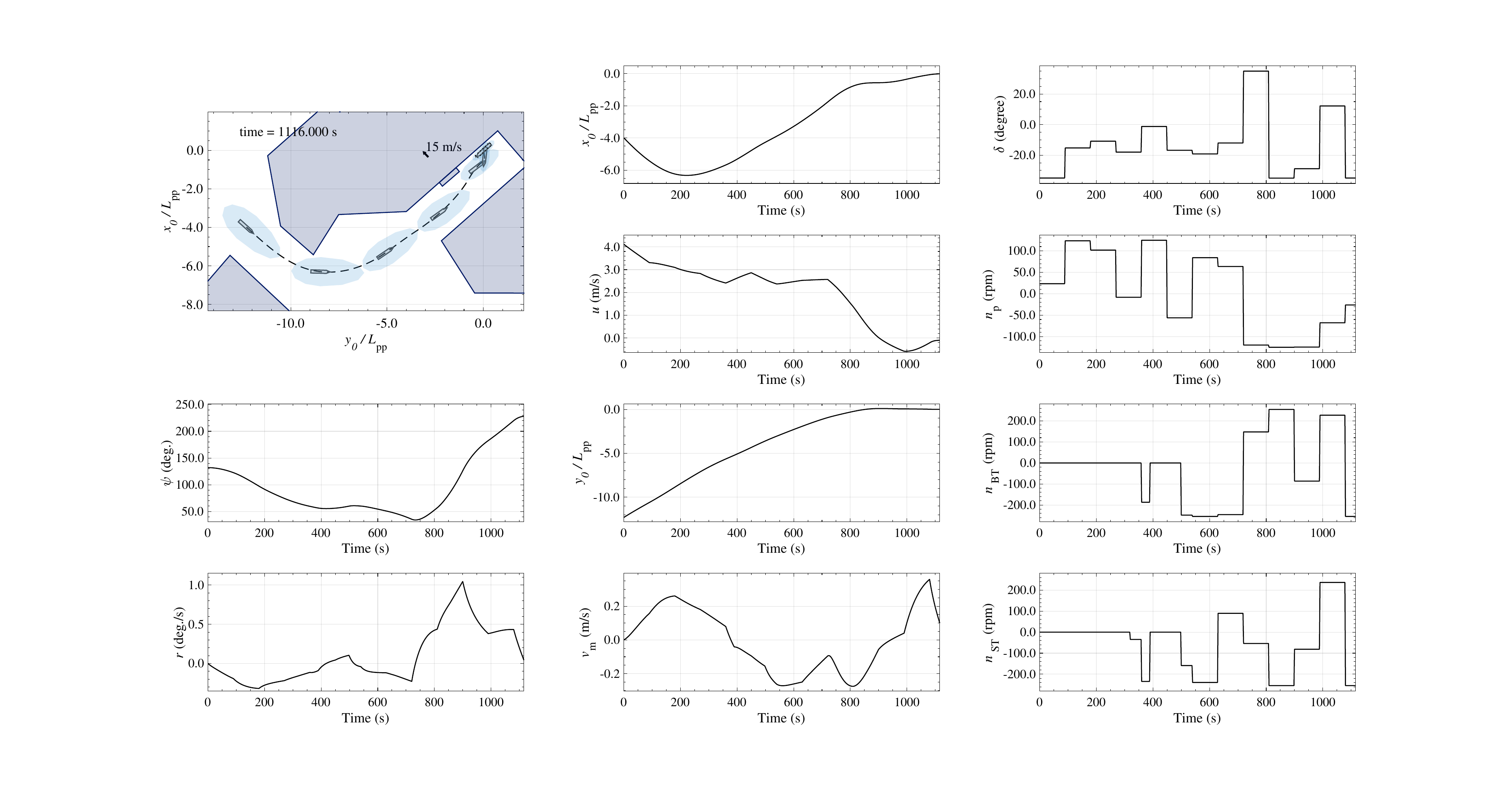}
    \caption{Optimal control inputs and achieved state of berthing on port of  \textit{Nanko}. \Add{Computation} with wind disturbance $U_{T}=15$ m/s. \Add{Ellipse shown in light blue on figure of trajectory at upper left is the ship domain used in collision avoidance algorithm, and arrow at upper left figure indicates wind direction.} }
    \label{fig:nnk_in_wind}
\end{figure*}

\begin{figure*}
    \centering
    \includegraphics[width=\linewidth]{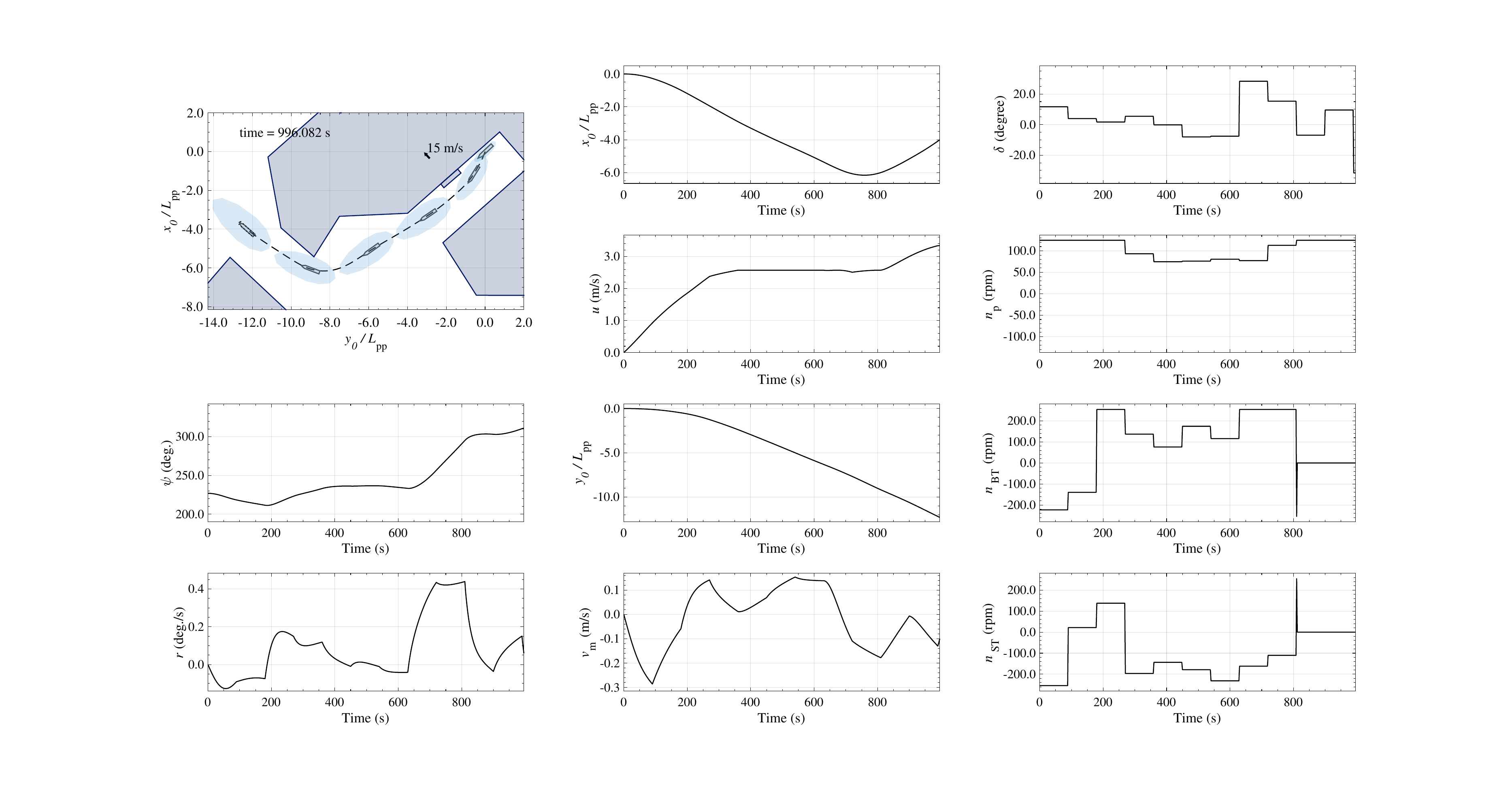}
    \caption{Optimal control inputs and achieved state of unberthing on port of  \textit{Nanko}. \Add{Computation} with wind disturbance $U_{T}=15$ m/s. \Add{Ellipse shown in light blue on figure of trajectory at upper left is the ship domain used in collision avoidance algorithm, and arrow at upper left figure indicates wind direction.} }
    \label{fig:nnk_out_wind}
\end{figure*}

Under the ideas discussed above, calculations with wind disturbance were conducted for both berthing and unberthing at  \textit{Nanko}. Fig.\ref{fig:nnk_in_wind}, \ref{fig:nnk_out_wind} show the berthing and unberthing results under wind disturbance of $\boldsymbol{\omega} = (137^{\circ},\ 15 \ \mathrm{m/s})$; the true wind direction $\gamma_{T}$ is perpendicular to the berth, which pushes the ship towards berth wall. From the figures, we can see that the CMA-ES can obtain a reasonable trajectory even under the wind disturbance when the actuator is capable enough.


\subsection{Verification on Different Port Geometry}\label{sec:ariake}
\begin{figure*}
    \centering
    \includegraphics[width=\linewidth]{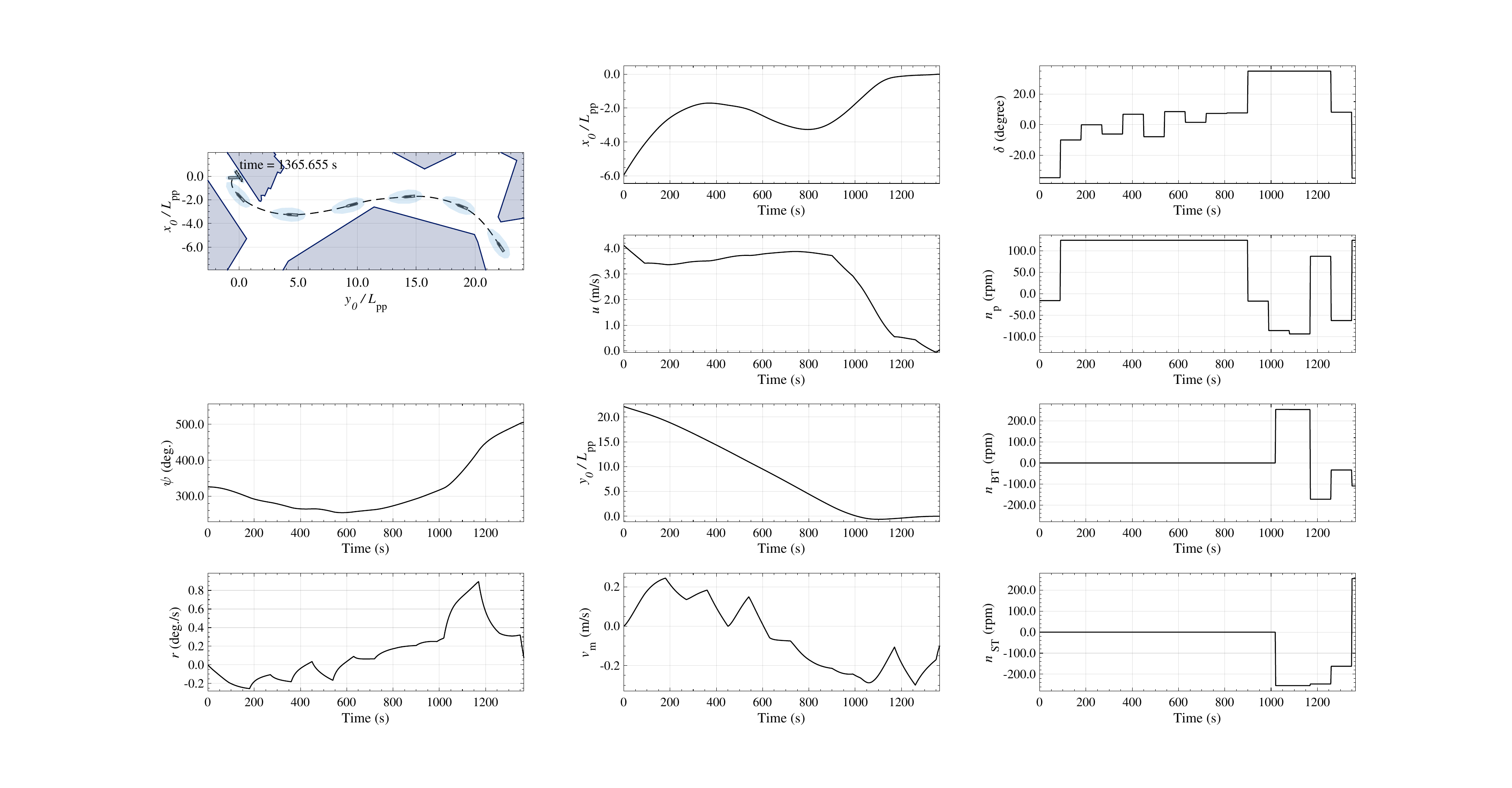}
    \caption{Optimal control inputs and achieved state  of berthing on port of  \textit{Ariake}. \Add{Computation} without the wind disturbance. \Add{Ellipse shown in light blue on figure of trajectory at upper left is the ship domain used in collision avoidance algorithm}}
    \label{fig:aak_in}
\end{figure*}

\begin{figure*}
    \centering
    \includegraphics[width=\linewidth]{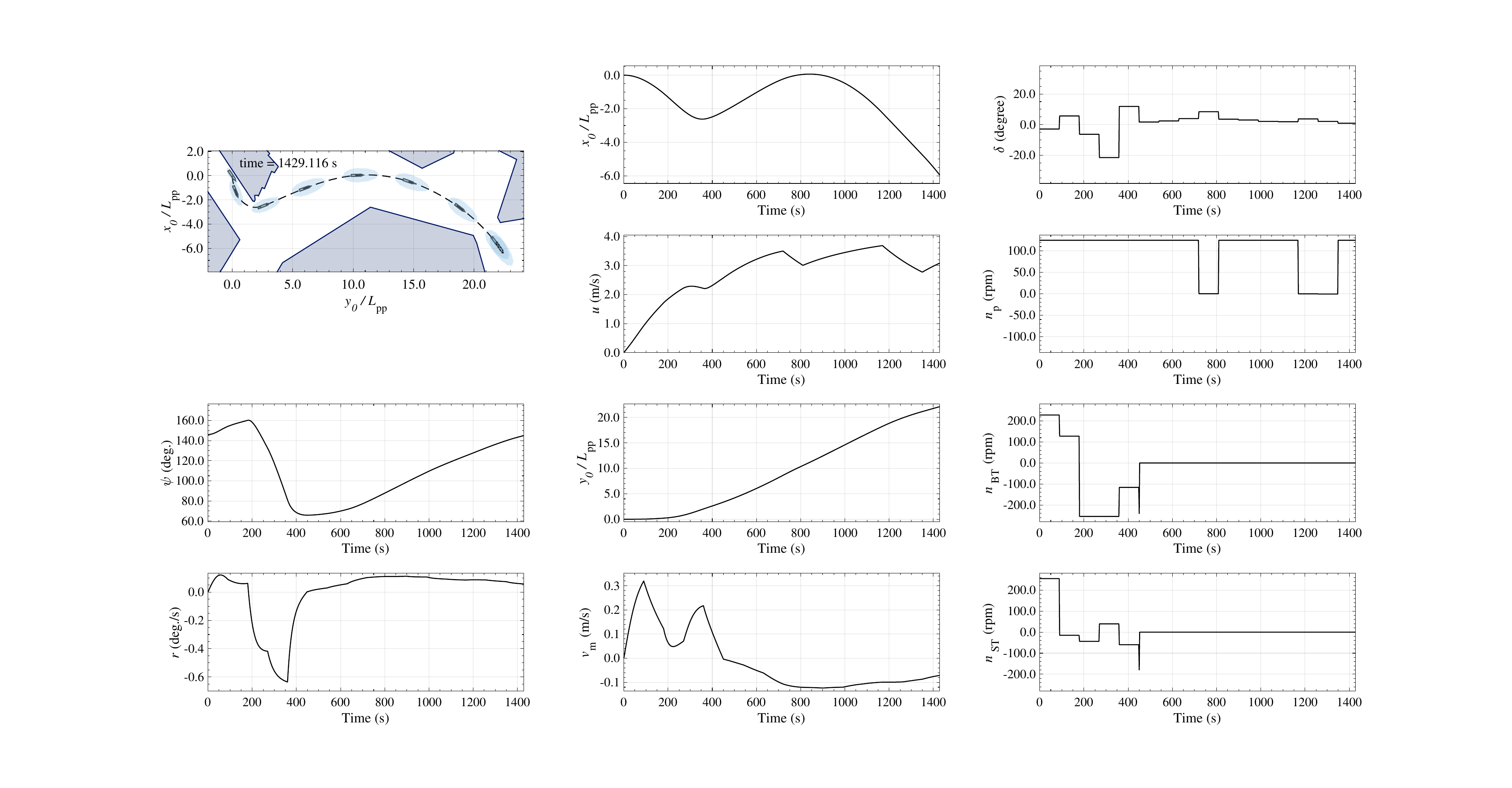}
    \caption{Optimal control inputs and achieved state  of unberthing on port of \textit{Ariake}. \Add{Computation} without the wind disturbance.}
    \label{fig:aak_out}
\end{figure*}

To show that the proposed method is applicable to multiple ports, it is necessary to verify the method on ports other than  \textit{Nanko}, used in the previous section.  \textit{Ariake} has the same spatial constraints characteristics while the travel distance is slightly longer than  \textit{Nanko} and berthing to the port side. Fig.\ref{fig:aak_in} and Fig.\ref{fig:aak_out} show the results, which seem reasonable. From the results, we can see that the proposed method is appreciable to different kinds of ports.

\subsection{Trajectory Planning Optimization with Waypoints}
\begin{figure*}
    \centering
    \includegraphics[width=\linewidth]{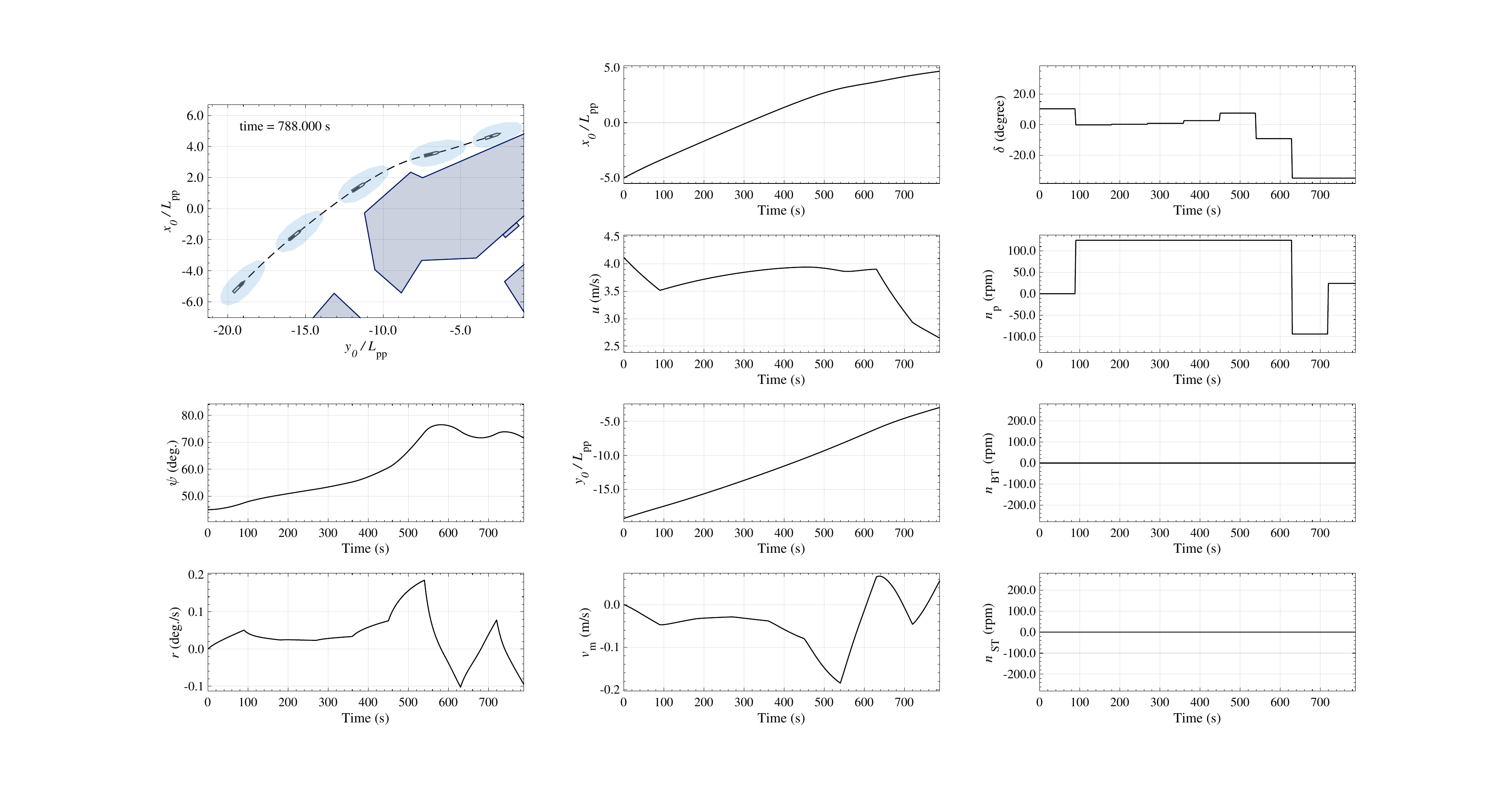}
    \caption{Example of unfavorable solution at  \textit{Nanko}. The trajectory did not enter the opening of breakwater to enter the berth.}
    \label{fig:nnk_withoutWP}
\end{figure*}

As shown in the previous section, the proposed method was able to find an appropriate trajectory even under the constraint of complex port geometry.
However, it was found that the proposed method does not provide a reasonable solution within the iteration limit when the ship starts from a position far from the breakwaters at the entrance of the port and the initial heading does not point to the entrance of the port. \cref{fig:nnk_withoutWP} shows one of the results of such cases, only the local minimum solution that the $t_{\mathrm{f}}$ ends early while advance speed $u$ remain large, was obtained. We assumed that this was because: solution exploration not proceed to solutions that navigate to the entrance of port, because the gradient of increasing $J$ occur around the solution toward the opening of the breakwater, due to the small choice of control inputs which can enter to the opening of the breakwater, and surrounding solutions of the favorable solution have significantly large $J$ caused by penalty of collision to breakwaters;  a solution that goes straight forward from the initial heading is likely to occur, even if the rudder angle is searched randomly. This is because to change the course of a ship, the rudder angle needs to maintain a certain period until the maneuver develops. 
% in these cases, the gradient of the present objective function shown on \cref{eq:objfunc} get steep when only the ship enters the narrow fairway between breakwaters which have a small choice of control inputs. Otherwise, the ship collides with surrounding breakwaters, making the objective function very large; hence, the solutions that not try to enter the openings will be a local minima.
% ・着桟地点へ向かうためには防波堤の開口部を衝突なく通過する必要があるが，そのための制御入力の組み合わせは限られている．
% ・また，防波堤の開口部を通過できる解の周辺解は障害物に衝突する解であるため，評価関数が非常に大きくなる．
% ・よって，終端状態量と衝突のみを関数とする本研究の評価関数では，岸壁の開口部に向かう解周辺の勾配が評価関数を増大させる方向であるので，開口部に進む制御入力の探索が進みにくい．
% ・加えて，船の方位が変わるまで運動が発達するまでには舵角をある程度連続して同方向に操舵し続ける必要があるので，初期値点で勾配がでない本研究の評価関数では，ランダムに舵角を変更したとしても行き足がついている初期状態では初期地点から直進する解が発生しやすい．
% ・上記2点が，Fig16で示した防波堤の開口部とは異なる方向に進む局所最適解が得られた原因と推察する．
% （次の段落で）
% ・このような岸壁の開口部と初期方位の条件のくみ合わせでも最適化を進めるためには，開口部に向けて進む解で評価関数が減少する勾配が生まれるようにする必要がある．
% ・そのために，評価関数にウェイポイントを追加する修正を加える．

To obtain reasonable solutions of trajectory planning  regardless of choice of the initial condition, we tried to generate favorable gradient of the objective function of \cref{eq:objfunc} when searching the solution around the entrance of port,  by including the waypoint portion to $J$. The $J$ shown on \cref{eq:objfunc} was modified to include the component of waypoints $J_{\text{WP}}$:
\begin{align}
\begin{aligned}
    \label{eq:objfunc_wp}
    J=\left(J_{1} +J_{\text{WP}}\right) \cdot t_{\mathrm{f}}+w_{c} C,
\end{aligned}
\end{align}
where
\begin{multline}\label{eq:objfunc_wp2}
    J_{1} = \sum_{i=1}^{6} w_{\mathrm{dim},i}\bigg(x_{\mathrm{tol},i}^2 \mathbf{1}_{\{\lvert x_{\mathrm{des}, i}-x_{i}(t_{\mathrm{f}})\rvert \leq x_{\mathrm{tol},i} \}} \\ 
    + w_\mathrm{pen} \big(x_{\mathrm{des}, i}-x_{i}(t_{\mathrm{f}})\big)^{2} \mathbf{1}_{\{\lvert x_{\mathrm{des}, i}-x_{i}(t_{\mathrm{f}})\rvert > x_{\mathrm{tol},i} \}} \bigg)\enspace,
\end{multline}
% \begin{align}
% \begin{aligned}
%     J_{1} &= \sum_{i}^{6} \max \left[ w_{\mathrm{pen}} \boldsymbol{w}_{\mathrm{dim}, i} \cdot\left\{\boldsymbol{x}_{\mathrm{des}, i}-\boldsymbol{x}_{i}(t_{\mathrm{f}})\right\}^{2}, \right.\\
%     & \left. \boldsymbol{x}_{\mathrm{dim}, i}\cdot \boldsymbol{x}_{\mathrm{tol}, i}^{2}\right], 
% \end{aligned}
% \end{align}
\begin{multline}\label{eq:objfunc_wp3}
    J_{\text{WP}} = \sum_{n_{\text{WP}}}\ x_{\mathrm{dim}, 1}\bigg( L_{\mathrm{tol}}^{2} 
    \mathbf{1}_{ \{L_{\mathrm{WP},i}\leq L_{\mathrm{tol}} \}} \\
    + L_{\mathrm{WP},i}^{2} \mathbf{1}_{ \{L_{\mathrm{WP},i} > L_{\mathrm{tol}} \}}  \bigg)\enspace,
\end{multline}
% \begin{align}
%     \label{eq:objfunc_wp3}
%     J_{\text{WP}} = \sum_{n_{\text{WP}}}\max\left\{\boldsymbol{x}_{\mathrm{dim}, 1}\cdot L_{\mathrm{WP},i}^{2}, \ \boldsymbol{x}_{\mathrm{dim}, 1}\cdot L_{\mathrm{tol}}^{2} \right\}, 
% \end{align}
\begin{equation}
     L_{\mathrm{WP},i} = \min\bigg(\sqrt{\{\boldsymbol{x}_{0}(t)-x_{\mathrm{WP},i}\}^2 + \{\boldsymbol{y}_{0}(t)-y_{\mathrm{WP},i}\}^2}\bigg)\enspace .
\end{equation}
Here, $n_{\text{WP}}$ is the total number of way point,  $x_{\mathrm{WP},i}, \ y_{\mathrm{WP},i}$ are the $x_{0}$ and $y_{0}$ coordinate of $i$-th waypoint, and $L_{\mathrm{WP},i}$ is the minimum distance between path of midship for $i$-th waypoint, respectively.
$L_{\mathrm{tol}}$ is the tolerance distance from way point, same as $\boldsymbol{x}_{\text{tol}}$. In this study, $L_{\mathrm{tol}}$ was set to $0.5 \ L_{\mathrm{pp}}$.

\begin{figure}
    \centering
    \includegraphics[width=\linewidth]{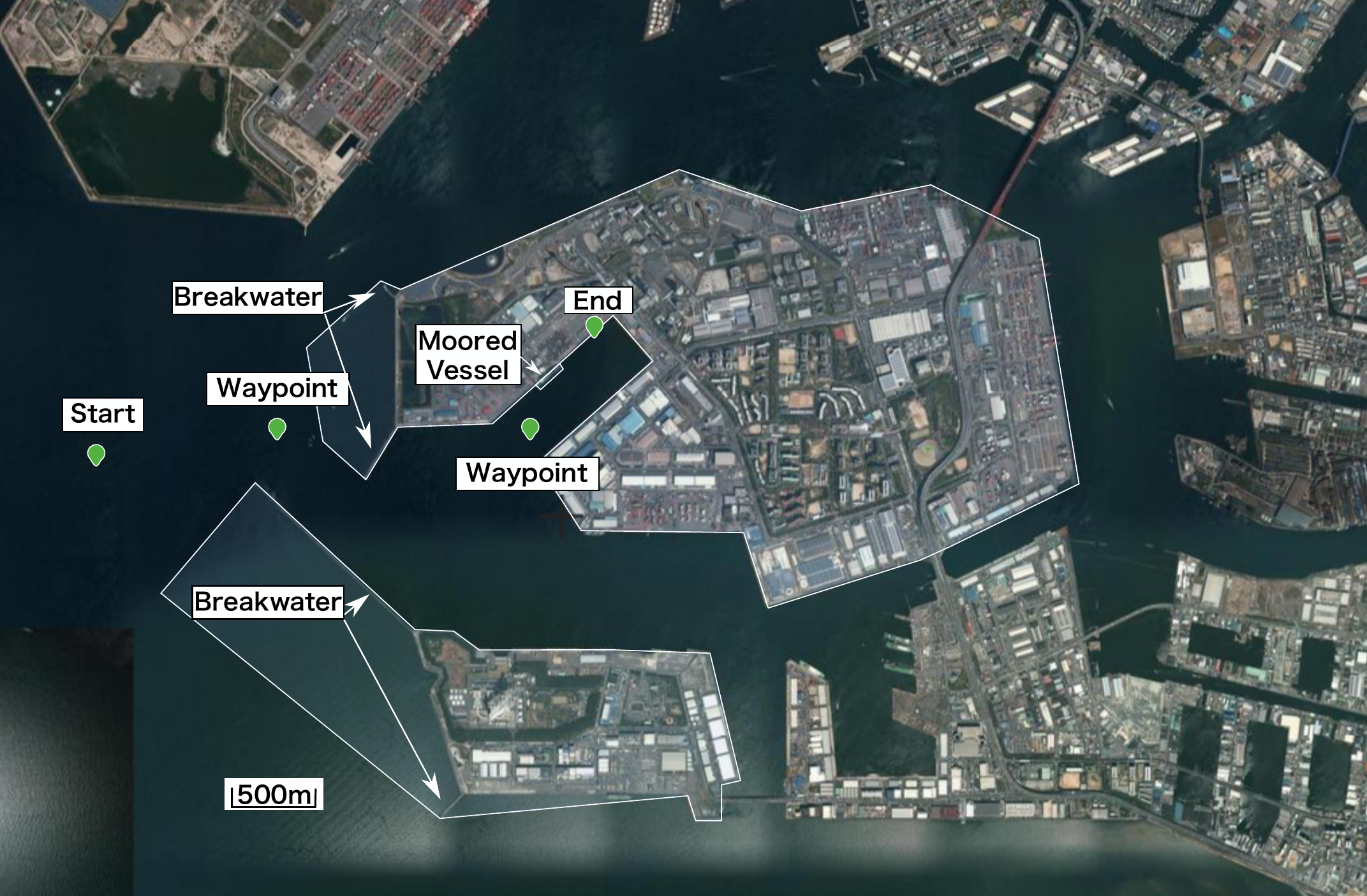}
    \caption{way point of port of  \textit{Nanko}}
    \label{fig:nnk_satimage_wp}
\end{figure}

\Add{On the optimization with waypoints and new objective function \Cref{eq:objfunc_wp}, we tested with two waypoints $n_{\text{WP}}=2$ on the berthing at the port of Nanko. Ideas of the selection of waypoints are as follows. First, waypoints are intended to generate a favorable gradient on $J$. We assumed that only a few waypoints were necessary to generate gradients, and designation of many waypoints undermines the automation of trajectory planning. Hence, the minimum number of waypoints ($n_{\text{WP}}=2$), which can prove that the proposed method can handle multiple points, was chosen. Second, the location of waypoints was chosen to provide a guide point of a reasonable path. Locations of waypoints are shown in \Cref{fig:nnk_satimage_wp}. The first waypoint was set at the harbor entrance, identical to the start point shown in \Cref{fig:Nanko_photo}. Another waypoint was set at the entrance of the confined berth. By setting  waypoints as a guide point of the reasonable path, CMA-ES tries to generate the solution which likely to pass nearby the waypoints. Hence it increases the likelihood of obtaining the preferable solution.}

The start point \Add{of the berthing with waypoints was} set at the passage outside the breakwater, shown on Fig.\ref{fig:nnk_satimage_wp}, \Add{which is identical to the computation can not obtain proper solution without waypoinst (\Cref{fig:nnk_withoutWP})}. Details of  initial condition are as follows:
\begin{equation}
    \label{eq:init_wp}
    \boldsymbol{x}_{\mathrm{int}} =  (-753.7 \ \mathrm{m}, 8.0 \ \mathrm{kn}, -2892.1 \ \mathrm{m}, 0.0 \ \mathrm{m/s}, 45 ^{\circ}, 0.0^{\circ}/\mathrm{s})^\mathsf{T} \enspace .
\end{equation}

\begin{figure*}
    \centering
    \includegraphics[width=\linewidth]{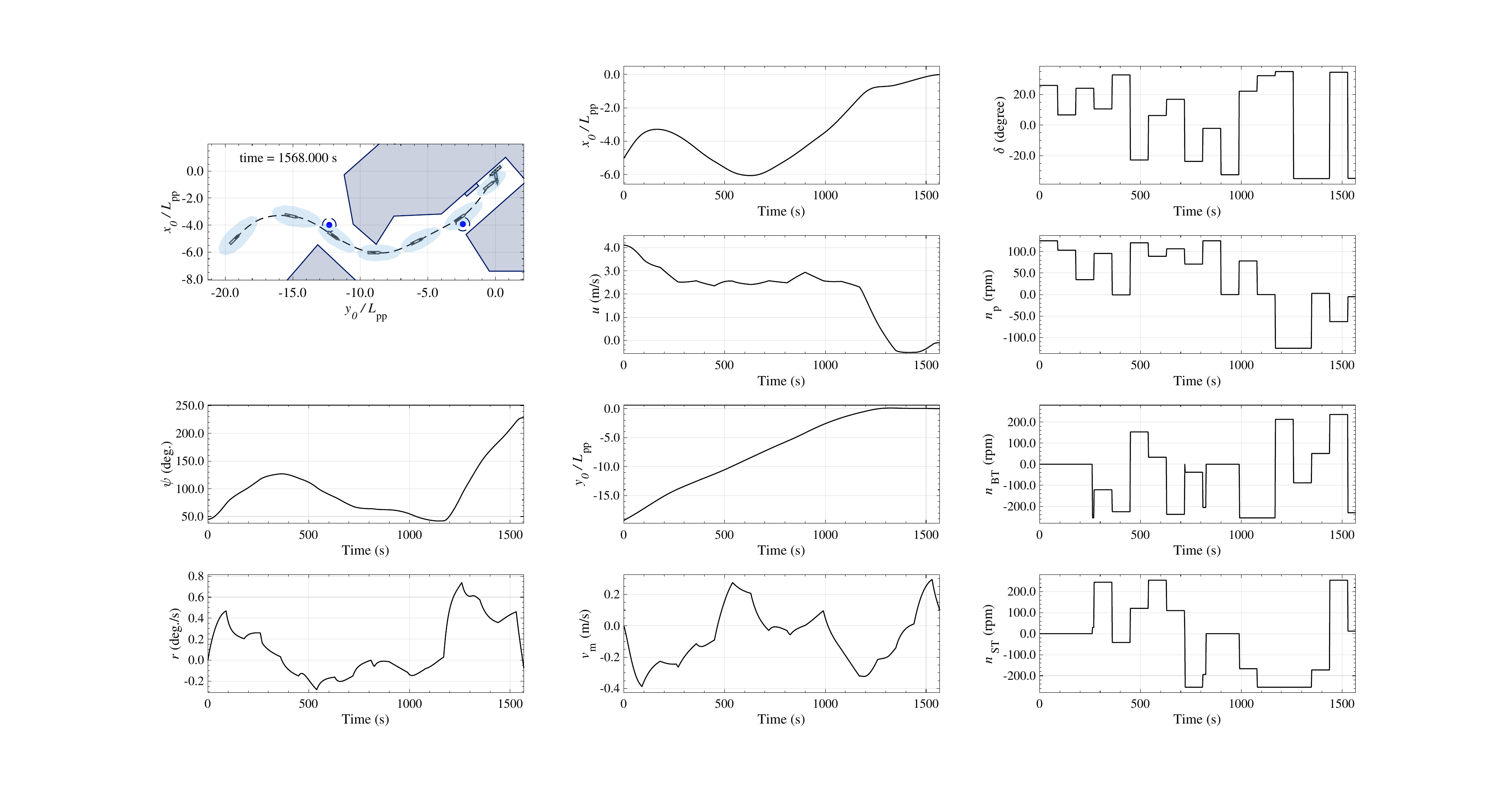}
    \caption{Computation with way points at  \textit{Nanko}. blue dashed circle shows the way point with tolerance distance.}
    \label{fig:nnk_withWP}
\end{figure*}

Computation results with waypoint are shown in Fig.\ref{fig:nnk_withWP}. By introducing waypoints and objective function with it, trajectory planning succeeded in which the ship arrived at the designated berthing point and properly tracing the waypoint. 


\section{Discussion}
% 手法の優位性をかく
As shown in the previous section, the proposed method can obtain collision-free, optimal trajectory and control input on real port geometry. \Add{Although the previous work of proposed method \citep{Maki2020,Maki2020b} was not able to incorporate complex spatial constraints and not considered the wind disturbance, the preset paper showed the improvement to overcome those limitations. We show applications become possible by the improvement later on in this section.} Our proposed method can directly apply to the arbitrary spatial constraints while similar research mentioned in  \cref{related} \citep{Martinsen2020, Bitar2020, Bergman2020} have adopted the two-stage method; discrete the continuous state space to apply graph search which obtains an initial guess. The present method can handle the comprehensive ship's dynamics \Add{as single-stage optimization} while graph search method can only address purely geometric space or limited way of ship dynamics. In addition, the present method varies the ship domain with speed which is natural for navigators' sense, while similar research \citep{Bitar2020, Bergman2020} have only adopted to circular safety region around the ship with a constant \Erase{radius}\Add{region size}.

% 得られたoffline trajectortの活用方法を書
Obtained trajectory and control can serve as a reference to the tracking control of berthing\Erase{although the present method can not directly apply to the berthing control because of the computational requirement of the optimization process, which is due to the iteration}. On the optimization of real-time control method with finite prediction horizon, such as model predictive control\citep{Li2020a}, those method could be unstable when the perdition horizon gets longer. However, by offering predefined trajectory as a reference, prediction horizon can be taken shorter. Trajectory obtained by the proposed method is suitable for reference because it is time-optimized collision-free trajectory.

\Add{Moreover, the proposed method} can also use to estimate the limit of actuator capability\Add{. By exploring the maximum wind velocity }which could maintain collision-free trajectory on certain \Add{actuator configuration} \Erase{wind} and port geometry, \Add{the ship designer can evaluate current design is sufficient to satisfy operation requirement or not}\Erase{ by testing on various wind conditions}. \Add{For example, results on \Cref{sec:nnk_wind} shows that subject ship's design can berth and unberth under the wind condition of $\boldsymbol{\omega} = (137^{\circ},\ 15~ \mathrm{m/s})$, while originally the maximum of side thruster's thrust was set to be equal to the wind pressure of $30~ \mathrm{m/s}$ lateral wind. However, the operational limit might be lower than $30~ \mathrm{m/s}$, because the ship has to maneuver inside the spatial constraint of port geometry. We will be able to find nominal operational limit by increasing the wind condition from $15~ \mathrm{m/s}$ and find the maximum of the present method can find the collision-free trajectory.}

\Add{The major drawback is its computational requirement. The proposed method can not directly apply to the real-time berthing control because of the computational requirement of the optimization process, which is due to the iteration. The proposed method takes few days for optimization. However, the proposed method can shorten its computational time requirement by applying parallel computation. CMA-ES uses multiple candidate solutions for a single iteration, as shown in \Cref{fig:flowchart}, so each candidate solution's computation on objective function can easily parallelize.  In addition, because of the computational requirement, we intended to use the proposed method combined with real-time control method as stated in the previous paragraph.}

\Add{Major remaining issue is that the proposed method is not suitable for collision avoidance with the passing ships. Theoretically, the proposed method can search for a collision-free trajectory if the trajectory of the dynamic obstacle was predetermined. However, this computation would  contribute little to autonomous berthing because on autonomous navigation, it is better to assume the passing ship's motion has some uncertainty. To overcome those drawbacks, The authors proposed using the current method to generate reference trajectory without passing ship, and the proposed method must be integrated with real-time control method for avoidance of passing ship. This issue is a future work for proposed method.}

\AddTwo{Another remaining issue of the proposed study is the treatment of unsteady wind disturbance. In this study, only steady wind disturbance was considered. As with the collision avoidance with passing ships, wind fluctuation must be overcome by integrating with the real-time control method.}

% 良いリファレンス経路であるための今後の課題をかく
Regarding the other potential future work, by improving two essential factors of trajectory planning, our work can add value as a reference trajectory; safety and optimization. From the safety perspective, the mathematical model of ship maneuver was assumed accurate enough. However, mathematical modeling of berthing maneuver itself still remains as research topic due to: the complexity and non-linearity of berthing motion; scale effect on hydrodynamic coefficients; and effect of shallow water and a sidewall. To negotiate the uncertainty of motion estimation caused by those factors, the robust optimization method is one option for improvement\AddTwo{, other than improving the mathematical model itself by utilizing Computational fluid dynamics or model experiments \citep{Ueno2015,Ueno2017,Chen2021}}. In robust optimization, the range of possible range of uncertainly on input data is set in advance, and optimization is performed on the input data, which gives the worst result \citep{ben2009robust}. From the optimization perspective, on the other hand, the practical berthing process is not just a single objective problem of time; optimization of berthing is a multi-objective optimization problem with the trade-off between time, energy consumption, and ride quality especially for passenger vessels. It is better to consider optimizing the berthing (and unberthing) trajectory as a multi-objective optimization problem that minimizes several factors simultaneously.

\section{Conclusion}
To achieve optimization of trajectory planning of berthing and unberthing of a ship at a real port, it is necessary to consider the spatial constraints, such as berths, breakwaters, buoys, or anchoring vessels. The novel collision avoidance algorithm was proposed, representing the spatial obstacles as polygons and includes a ship domain to maintain the sufficient distance to obstacles. The proposed ship domain consists of domain boundary vertices that discretize the ellipse shape of the ship domain with a uniform argument angle, and the size of the ship domain changes with ship's speed, which is based on the navigators' knowledge and experience.

Optimization of trajectory planning for two existing ports,  \textit{Nanko} and  \textit{Ariake}, were conducted by using the proposed collision avoidance algorithm with CMA-ES and the MMG model. By introducing the proposed collision avoidance algorithm, collision-free time-optimal control inputs and trajectories were obtained. The effect of wind disturbance was also considered. By considering the wind force, the obtained trajectory will be more robust that considers the limits of actuators. Additionally, optimization with waypoints also shown in this study. Including the waypoints to objective function makes it more robust to the choice of the initial location of berthing simulation. Although its long computation time due to the iteration, it should be noted that the present method requires only the mathematical model of ship maneuver, a start point, an end point, polygons of spatial constraints, and waypoints, if necessary, to obtain an optimal trajectory at the real port. The optimal trajectory can serve as a predefined reference to the real-time tracking control of berthing and unberthing.

\printcredits

\section*{Declaration of competing interest}
The authors declare that they have no known competing financial
interests or personal relationships that could have appeared to influence the work reported in this paper.

\section*{Acknowledgment}
This study was supported by a Grant-in-Aid for Scientific Research from the Japan Society for Promotion of Science (JSPS KAKENHI Grant \#19K04858).
The study also received assistance from JFY2018 Fundamental Research Developing Association for Shipbuilding and Offshore (REDAS) in Japan.

This is accepted manuscript. Published journal Article is shown on \url{https://doi.org/10.1016/j.oceaneng.2021.110390}.
\appendix

\section{Supplementary data}
% \textcolor{blue}{Comment: The author would like to add video files as supplementary data on this chapter.}

% (text for electric version)
% The following is the supplementary data related to this article:

% \begin{itemize}
%     \item Video 1: Animation of optimal control input and trajectory of berthing on the port of  \textit{Nanko}. Without wind disturbance \AddTwo{and with a comparison computation of previous method}.
%     \item Video 2: Animation of optimal control input and trajectory of unberthing on the port of  \textit{Nanko}. Without wind disturbance.
%     \item Video 3: Animation of optimal control input and trajectory of berthing on the port of  \textit{Nanko}. With wind disturbance $U_{T}=15$ m/s.
%     \item Video 4: Animation of optimal control input and trajectory of unberthing on the port of  \textit{Nanko}. With wind disturbance $U_{T}=15$ m/s.
%     \item Video 5: Animation of optimal control input and trajectory of berthing on the port of  \textit{Ariake}. Without wind disturbance.
%     \item Video 6: Animation of optimal control input and trajectory of unberthing on the port of  \textit{Ariake}. Without wind disturbance.
%     \item Video 7: Animation of optimal control input and trajectory of berthing on the port of  \textit{Nanko} with waypoints.
% \end{itemize}
% (text for the print version )
Supplementary data related to this article can be found at \url{https://doi.org/10.1016/j.oceaneng.2021.110390}.
\appendix


\begin{thebibliography}{44}
    \expandafter\ifx\csname natexlab\endcsname\relax\def\natexlab#1{#1}\fi
    \providecommand{\url}[1]{\texttt{#1}}
    \providecommand{\href}[2]{#2}
    \providecommand{\path}[1]{#1}
    \providecommand{\DOIprefix}{doi:}
    \providecommand{\ArXivprefix}{arXiv:}
    \providecommand{\URLprefix}{URL: }
    \providecommand{\Pubmedprefix}{pmid:}
    \providecommand{\doi}[1]{\href{http://dx.doi.org/#1}{\path{#1}}}
    \providecommand{\Pubmed}[1]{\href{pmid:#1}{\path{#1}}}
    \providecommand{\bibinfo}[2]{#2}
    \ifx\xfnm\relax \def\xfnm[#1]{\unskip,\space#1}\fi
    %Type = Phdthesis
    \bibitem[{Ahmed(2015)}]{Ahmed2015}
    \bibinfo{author}{Ahmed, Y.A.}, \bibinfo{year}{2015}.
    \newblock \bibinfo{title}{Automatic Berthing Control Practically Applicable
      under Wind Disturbances}.
    \newblock Ph.D. thesis. Osaka University.
    %Type = Article
    \bibitem[{Amendola et~al.(2020)Amendola, Miura, Costa, Cozman and
      Tannuri}]{Amendola2020}
    \bibinfo{author}{Amendola, J.}, \bibinfo{author}{Miura, L.S.},
      \bibinfo{author}{Costa, A.H.}, \bibinfo{author}{Cozman, F.G.},
      \bibinfo{author}{Tannuri, E.A.}, \bibinfo{year}{2020}.
    \newblock \bibinfo{title}{Navigation in restricted channels under environmental
      conditions: Fast-time simulation by asynchronous deep reinforcement
      learning}.
    \newblock \bibinfo{journal}{IEEE Access} \bibinfo{volume}{8},
      \bibinfo{pages}{149199--149213}.
    \newblock \DOIprefix\doi{10.1109/ACCESS.2020.3015661}.
    %Type = Inproceedings
    \bibitem[{Auger and Hansen(2005)}]{Auger2005}
    \bibinfo{author}{Auger, A.}, \bibinfo{author}{Hansen, N.},
      \bibinfo{year}{2005}.
    \newblock \bibinfo{title}{A restart cma evolution strategy with increasing
      population size}, in: \bibinfo{booktitle}{2005 IEEE Congress on Evolutionary
      Computation}, \bibinfo{publisher}{IEEE}. pp. \bibinfo{pages}{1769--1776}.
    \newblock \URLprefix \url{http://ieeexplore.ieee.org/document/1554902/},
      \DOIprefix\doi{10.1109/CEC.2005.1554902}.
    %Type = Book
    \bibitem[{Ben-Tal et~al.(2009)Ben-Tal, Ghaoui and Nemirovski}]{ben2009robust}
    \bibinfo{author}{Ben-Tal, A.}, \bibinfo{author}{Ghaoui, L.},
      \bibinfo{author}{Nemirovski, A.}, \bibinfo{year}{2009}.
    \newblock \bibinfo{title}{Robust Optimization}.
    \newblock Princeton Series in Applied Mathematics,
      \bibinfo{publisher}{Princeton University Press}.
    \newblock \URLprefix \url{https://books.google.co.jp/books?id=DttjR7IpjUEC}.
    %Type = Inproceedings
    \bibitem[{Bergman et~al.(2020)Bergman, Ljungqvist, Linder and
      Axehill}]{Bergman2020}
    \bibinfo{author}{Bergman, K.}, \bibinfo{author}{Ljungqvist, O.},
      \bibinfo{author}{Linder, J.}, \bibinfo{author}{Axehill, D.},
      \bibinfo{year}{2020}.
    \newblock \bibinfo{title}{{An Optimization-Based Motion Planner for Autonomous
      Maneuvering of Marine Vessels in Complex Environments}}, in:
      \bibinfo{booktitle}{2020 59th IEEE Conference on Decision and Control (CDC)},
      \bibinfo{publisher}{IEEE}. pp. \bibinfo{pages}{5283--5290}.
    \newblock \URLprefix \url{https://ieeexplore.ieee.org/document/9303746/
      http://arxiv.org/abs/2012.12145},
      \DOIprefix\doi{10.1109/CDC42340.2020.9303746},
      \href{http://arxiv.org/abs/2005.02674}{\tt arXiv:2005.02674}.
    %Type = Article
    \bibitem[{Bitar et~al.(2020)Bitar, Martinsen, Lekkas and Breivik}]{Bitar2020}
    \bibinfo{author}{Bitar, G.}, \bibinfo{author}{Martinsen, A.B.},
      \bibinfo{author}{Lekkas, A.M.}, \bibinfo{author}{Breivik, M.},
      \bibinfo{year}{2020}.
    \newblock \bibinfo{title}{Two-stage optimized trajectory planning for asvs
      under polygonal obstacle constraints: Theory and experiments}.
    \newblock \bibinfo{journal}{IEEE Access} \bibinfo{volume}{8},
      \bibinfo{pages}{199953--199969}.
    \newblock \URLprefix
      \url{https://www.maritime-executive.com/article/rolls-royce-and-wartsila-},
      \DOIprefix\doi{10.1109/access.2020.3035256}.
    %Type = Article
    \bibitem[{Chen et~al.(2021)Chen, Tu, Ti, Wang and Hu}]{Chen2021}
    \bibinfo{author}{Chen, G.}, \bibinfo{author}{Tu, J.}, \bibinfo{author}{Ti, X.},
      \bibinfo{author}{Wang, Z.}, \bibinfo{author}{Hu, H.}, \bibinfo{year}{2021}.
    \newblock \bibinfo{title}{Hydrodynamic model of the beaver-like bendable webbed
      foot and paddling characteristics under different flow velocities}.
    \newblock \bibinfo{journal}{Ocean Engineering} \bibinfo{volume}{234},
      \bibinfo{pages}{109179}.
    \newblock \URLprefix
      \url{https://www.sciencedirect.com/science/article/pii/S0029801821006120},
      \DOIprefix\doi{https://doi.org/10.1016/j.oceaneng.2021.109179}.
    %Type = Article
    \bibitem[{Fujii and Tanaka(1971)}]{fujii1971}
    \bibinfo{author}{Fujii, Y.}, \bibinfo{author}{Tanaka, K.},
      \bibinfo{year}{1971}.
    \newblock \bibinfo{title}{Traffic capacity}.
    \newblock \bibinfo{journal}{Journal of Navigation} \bibinfo{volume}{24},
      \bibinfo{pages}{543--552}.
    \newblock \URLprefix
      \url{https://www.cambridge.org/core/product/identifier/S0373463300022384/type/journal_article},
      \DOIprefix\doi{10.1017/S0373463300022384}.
    %Type = Article
    \bibitem[{Fujiwara et~al.(1998)Fujiwara, Ueno and Nimura}]{Fujiwara1998}
    \bibinfo{author}{Fujiwara, T.}, \bibinfo{author}{Ueno, M.},
      \bibinfo{author}{Nimura, T.}, \bibinfo{year}{1998}.
    \newblock \bibinfo{title}{Estimation of wind forces and moments acting on
      ships}.
    \newblock \bibinfo{journal}{Journal of the Society of Naval Architects of
      Japan} \bibinfo{volume}{1998}, \bibinfo{pages}{77--90}.
    \newblock \URLprefix
      \url{http://joi.jlc.jst.go.jp/JST.Journalarchive/jjasnaoe1968/1998.77?from=CrossRef},
      \DOIprefix\doi{10.2534/jjasnaoe1968.1998.77}.
    %Type = Article
    \bibitem[{Goodwin(1975)}]{goodwin1975}
    \bibinfo{author}{Goodwin, E.M.}, \bibinfo{year}{1975}.
    \newblock \bibinfo{title}{A statistical study of ship domains}.
    \newblock \bibinfo{journal}{Journal of Navigation} \bibinfo{volume}{28},
      \bibinfo{pages}{328–344}.
    \newblock \DOIprefix\doi{10.1017/S0373463300041230}.
    %Type = Masterthesis
    \bibitem[{Hachii(2004)}]{Hachii2004}
    \bibinfo{author}{Hachii, T.}, \bibinfo{year}{2004}.
    \newblock \bibinfo{title}{The prediciton of manoeuvring motion on ships with
      low speed using standard MMG model}.
    \newblock Master's thesis. Osaka University.
    %Type = Article
    \bibitem[{Hansen et~al.(2013)Hansen, Jensen, Lehn-Schiøler, Melchild,
      Rasmussen and Ennemark}]{Hansen2013}
    \bibinfo{author}{Hansen, M.G.}, \bibinfo{author}{Jensen, T.K.},
      \bibinfo{author}{Lehn-Schiøler, T.}, \bibinfo{author}{Melchild, K.},
      \bibinfo{author}{Rasmussen, F.M.}, \bibinfo{author}{Ennemark, F.},
      \bibinfo{year}{2013}.
    \newblock \bibinfo{title}{Empirical ship domain based on ais data}.
    \newblock \bibinfo{journal}{Journal of Navigation} \bibinfo{volume}{66},
      \bibinfo{pages}{931--940}.
    \newblock \URLprefix
      \url{https://www.cambridge.org/core/product/identifier/S0373463313000489/type/journal_article},
      \DOIprefix\doi{10.1017/S0373463313000489}.
    %Type = Inbook
    \bibitem[{Hansen(2006)}]{Hansen2006}
    \bibinfo{author}{Hansen, N.}, \bibinfo{year}{2006}.
    \newblock \bibinfo{title}{The CMA Evolution Strategy: A Comparing Review}.
      \bibinfo{publisher}{Springer Berlin Heidelberg}, \bibinfo{address}{Berlin,
      Heidelberg}.
    \newblock pp. \bibinfo{pages}{75--102}.
    \newblock \URLprefix \url{https://doi.org/10.1007/3-540-32494-1_4},
      \DOIprefix\doi{10.1007/3-540-32494-1_4}.
    %Type = Inproceedings
    \bibitem[{Hasegawa and Kitera(1993)}]{Hasegawa1993}
    \bibinfo{author}{Hasegawa, K.}, \bibinfo{author}{Kitera, K.},
      \bibinfo{year}{1993}.
    \newblock \bibinfo{title}{Mathematical model of manoeuvrability at low advance
      speed and its application to berthing control}, in: \bibinfo{booktitle}{2nd
      Japan-Korea Joint Workshop on Ship and Marine Hydrodynamics}, pp.
      \bibinfo{pages}{311--321}.
    %Type = Article
    \bibitem[{Inoue et~al.(1994)Inoue, Usami and Shibata}]{INOUE1994}
    \bibinfo{author}{Inoue, K.}, \bibinfo{author}{Usami, S.},
      \bibinfo{author}{Shibata, T.}, \bibinfo{year}{1994}.
    \newblock \bibinfo{title}{Modelling of mariners' senses on minimum passing
      distance between ships in harbour}.
    \newblock \bibinfo{journal}{The Journal of Japan Institute of Navigation}
      \bibinfo{volume}{90}, \bibinfo{pages}{297--306}.
    \newblock \DOIprefix\doi{10.9749/jin.90.297}.
    %Type = Article
    \bibitem[{Jensen et~al.(2013)Jensen, Hansen, Lehn-Schiøler, Melchild,
      Rasmussen and Ennemark}]{Jensen2013}
    \bibinfo{author}{Jensen, T.K.}, \bibinfo{author}{Hansen, M.G.},
      \bibinfo{author}{Lehn-Schiøler, T.}, \bibinfo{author}{Melchild, K.},
      \bibinfo{author}{Rasmussen, F.M.}, \bibinfo{author}{Ennemark, F.},
      \bibinfo{year}{2013}.
    \newblock \bibinfo{title}{Free flow–efficiency of a one-way traffic lane
      between two pylons}.
    \newblock \bibinfo{journal}{Journal of Navigation} \bibinfo{volume}{66},
      \bibinfo{pages}{941--951}.
    \newblock \URLprefix \url{https://doi.org/10.1017/S0373463313000362
      https://www.cambridge.org/core/product/identifier/S0373463313000362/type/journal_article},
      \DOIprefix\doi{10.1017/S0373463313000362}.
    %Type = Article
    \bibitem[{Kitagawa et~al.(2015)Kitagawa, Tsukada and Miyazaki}]{Kitagawa2015}
    \bibinfo{author}{Kitagawa, Y.}, \bibinfo{author}{Tsukada, Y.},
      \bibinfo{author}{Miyazaki, H.}, \bibinfo{year}{2015}.
    \newblock \bibinfo{title}{2015s-gs1-1 a study on mathematical models of
      propeller and rudder under maneuvering with propeller reverse rotation}.
    \newblock \bibinfo{journal}{Conference Proceedings The Japan Society of Naval
      Architects and Ocean Engineers} \bibinfo{volume}{20},
      \bibinfo{pages}{117--120}.
    \newblock \DOIprefix\doi{10.14856/conf.20.0_117}.
    %Type = Misc
    \bibitem[{Kobayashi(1988)}]{Kobayashi1988}
    \bibinfo{author}{Kobayashi, E.}, \bibinfo{year}{1988}.
    \newblock \bibinfo{title}{A simulation study on ship manoeuvrability at low
      speeds}.
    \newblock \bibinfo{howpublished}{Akishima Laboratory, Ocean Engineering
      Research Section, Mitsubishi Heave Industries Ltd. Published in: Mitsubishi
      Technical Bulletin No. 180}.
    %Type = Article
    \bibitem[{Kobayashi(1990)}]{Kobayashi1990}
    \bibinfo{author}{Kobayashi, E.}, \bibinfo{year}{1990}.
    \newblock \bibinfo{title}{Manoeuvring simulation at low speed for a ship with
      twin screw, rudders and thrusters}.
    \newblock \bibinfo{journal}{MARSIM and ICSM 90, Intl. Conference, Marine
      Simulation and Ship Manoeuvrability} .
    %Type = Article
    \bibitem[{Kobayashi et~al.(2002)Kobayashi, Blok, Barr, Kim and
      Nowicki}]{Kobayashi2002}
    \bibinfo{author}{Kobayashi, H.}, \bibinfo{author}{Blok, J.},
      \bibinfo{author}{Barr, R.}, \bibinfo{author}{Kim, Y.S.},
      \bibinfo{author}{Nowicki, J.}, \bibinfo{year}{2002}.
    \newblock \bibinfo{title}{The specialist committee on esso osaka final report
      and recommendations to the 23rd ittc}.
    \newblock \bibinfo{journal}{23rd International Towing Tank Conference}
      \bibinfo{volume}{II}, \bibinfo{pages}{581--743}.
    %Type = Article
    \bibitem[{Li et~al.(2020)Li, Liu, Negenborn and Wu}]{Li2020a}
    \bibinfo{author}{Li, S.}, \bibinfo{author}{Liu, J.},
      \bibinfo{author}{Negenborn, R.R.}, \bibinfo{author}{Wu, Q.},
      \bibinfo{year}{2020}.
    \newblock \bibinfo{title}{{Automatic Docking for Underactuated Ships Based on
      Multi-Objective Nonlinear Model Predictive Control}}.
    \newblock \bibinfo{journal}{IEEE Access} \bibinfo{volume}{8},
      \bibinfo{pages}{70044--70057}.
    \newblock \DOIprefix\doi{10.1109/ACCESS.2020.2984812}.
    %Type = Article
    \bibitem[{Liu et~al.(2016)Liu, Zhou, Li, Wang and Liu}]{Liu2016}
    \bibinfo{author}{Liu, J.}, \bibinfo{author}{Zhou, F.}, \bibinfo{author}{Li,
      Z.}, \bibinfo{author}{Wang, M.}, \bibinfo{author}{Liu, R.W.},
      \bibinfo{year}{2016}.
    \newblock \bibinfo{title}{{Dynamic Ship Domain Models for Capacity Analysis of
      Restricted Water Channels}}.
    \newblock \bibinfo{journal}{Journal of Navigation} \bibinfo{volume}{69},
      \bibinfo{pages}{481--503}.
    \newblock \URLprefix
      \url{https://www.cambridge.org/core/product/identifier/S0373463315000764/type/journal{\_}article},
      \DOIprefix\doi{10.1017/S0373463315000764}.
    %Type = Article
    \bibitem[{Maki et~al.(2020a)Maki, Akimoto and Naoya}]{Maki2020}
    \bibinfo{author}{Maki, A.}, \bibinfo{author}{Akimoto, Y.},
      \bibinfo{author}{Naoya, U.}, \bibinfo{year}{2020}a.
    \newblock \bibinfo{title}{Application of optimal control theory based on the
      evolution strategy (cma-es) to automatic berthing (part: 2)}.
    \newblock \bibinfo{journal}{Journal of Marine Science and Technology (Japan)}
      \DOIprefix\doi{10.1007/s00773-020-00774-x}.
    %Type = Article
    \bibitem[{Maki et~al.(2020b)Maki, Sakamoto, Akimoto, Nishikawa and
      Umeda}]{Maki2020b}
    \bibinfo{author}{Maki, A.}, \bibinfo{author}{Sakamoto, N.},
      \bibinfo{author}{Akimoto, Y.}, \bibinfo{author}{Nishikawa, H.},
      \bibinfo{author}{Umeda, N.}, \bibinfo{year}{2020}b.
    \newblock \bibinfo{title}{Application of optimal control theory based on the
      evolution strategy (cma-es) to automatic berthing}.
    \newblock \bibinfo{journal}{Journal of Marine Science and Technology}
      \bibinfo{volume}{25}, \bibinfo{pages}{221--233}.
    \newblock \URLprefix \url{https://doi.org/10.1007/s00773-019-00642-3
      http://link.springer.com/10.1007/s00773-019-00642-3},
      \DOIprefix\doi{10.1007/s00773-019-00642-3}.
    %Type = Article
    \bibitem[{Martinsen et~al.(2019)Martinsen, Lekkas and Gros}]{Martinsen2019}
    \bibinfo{author}{Martinsen, A.B.}, \bibinfo{author}{Lekkas, A.M.},
      \bibinfo{author}{Gros, S.}, \bibinfo{year}{2019}.
    \newblock \bibinfo{title}{Autonomous docking using direct optimal control}.
    \newblock \bibinfo{journal}{IFAC-PapersOnLine} \bibinfo{volume}{52},
      \bibinfo{pages}{97--102}.
    \newblock \DOIprefix\doi{10.1016/j.ifacol.2019.12.290}.
    %Type = Article
    \bibitem[{Martinsen et~al.(2021)Martinsen, Lekkas and Gros}]{Martinsen2020}
    \bibinfo{author}{Martinsen, A.B.}, \bibinfo{author}{Lekkas, A.M.},
      \bibinfo{author}{Gros, S.}, \bibinfo{year}{2021}.
    \newblock \bibinfo{title}{Optimal model-based trajectory planning with static
      polygonal constraints}.
    \newblock \bibinfo{journal}{IEEE Transactions on Control Systems Technology} ,
      \bibinfo{pages}{1--12}\DOIprefix\doi{10.1109/TCST.2021.3094617}.
    %Type = Article
    \bibitem[{Miyauchi et~al.(2021)Miyauchi, Maki, Umeda, Rachman and
      Akimoto}]{Miyauchi20201}
    \bibinfo{author}{Miyauchi, Y.}, \bibinfo{author}{Maki, A.},
      \bibinfo{author}{Umeda, N.}, \bibinfo{author}{Rachman, D.M.},
      \bibinfo{author}{Akimoto, Y.}, \bibinfo{year}{2021}.
    \newblock \bibinfo{title}{System parameter exploration of ship maneuvering
      model for automatic docking / berthing using cma-es}.
    \newblock \bibinfo{journal}{arXiv preprint},
    \newblock \URLprefix
      \url{https://arxiv.org/abs/2111.06124}.
    %Type = Techreport
    \bibitem[{Murakami et~al.(2015)Murakami, Takenobu, Miyata and
      Yoneyama}]{Murakami2015}
    \bibinfo{author}{Murakami, K.}, \bibinfo{author}{Takenobu, M.},
      \bibinfo{author}{Miyata, T.}, \bibinfo{author}{Yoneyama, H.},
      \bibinfo{year}{2015}.
    \newblock \bibinfo{title}{Fundamental Analysis on the Characteristics of
      Berthing Velocity of Ships For the Design of Port Facilities}.
    \newblock \bibinfo{type}{Technical Report}. TECHNICAL NOTE of National
      Institute for Land and Infrastructure Management.
    %Type = Article
    \bibitem[{Ogawa and Kasai(1978)}]{Ogawa1978}
    \bibinfo{author}{Ogawa, A.}, \bibinfo{author}{Kasai, H.}, \bibinfo{year}{1978}.
    \newblock \bibinfo{title}{{On the mathematical model of manoeuvring motion of
      ships}}.
    \newblock \bibinfo{journal}{International Shipbuilding Progress}
      \bibinfo{volume}{25}, \bibinfo{pages}{306--319}.
    \newblock \DOIprefix\doi{10.3233/ISP-1978-2529202}.
    %Type = Article
    \bibitem[{Pietrzykowski(2008)}]{Pietrzykowski2008}
    \bibinfo{author}{Pietrzykowski, Z.}, \bibinfo{year}{2008}.
    \newblock \bibinfo{title}{Ship's fuzzy domain – a criterion for navigational
      safety in narrow fairways}.
    \newblock \bibinfo{journal}{Journal of Navigation} \bibinfo{volume}{61},
      \bibinfo{pages}{499--514}.
    \newblock \URLprefix
      \url{https://www.cambridge.org/core/product/identifier/S0373463308004682/type/journal_article},
      \DOIprefix\doi{10.1017/S0373463308004682}.
    %Type = Article
    \bibitem[{Roubos et~al.(2017)Roubos, Groenewegen and Peters}]{Roubos2017}
    \bibinfo{author}{Roubos, A.}, \bibinfo{author}{Groenewegen, L.},
      \bibinfo{author}{Peters, D.J.}, \bibinfo{year}{2017}.
    \newblock \bibinfo{title}{Berthing velocity of large seagoing vessels in the
      port of rotterdam}.
    \newblock \bibinfo{journal}{Marine Structures} \bibinfo{volume}{51},
      \bibinfo{pages}{202--219}.
    \newblock \DOIprefix\doi{10.1016/j.marstruc.2016.10.011}.
    %Type = Inproceedings
    \bibitem[{Sakamoto and Akimoto(2017)}]{Sakamoto2017}
    \bibinfo{author}{Sakamoto, N.}, \bibinfo{author}{Akimoto, Y.},
      \bibinfo{year}{2017}.
    \newblock \bibinfo{title}{Modified box constraint handling for the covariance
      matrix adaptation evolution strategy}, in: \bibinfo{booktitle}{Proceedings of
      the Genetic and Evolutionary Computation Conference Companion},
      \bibinfo{publisher}{Association for Computing Machinery},
      \bibinfo{address}{New York, NY, USA}. p. \bibinfo{pages}{183–184}.
    \newblock \URLprefix \url{https://doi.org/10.1145/3067695.3075986},
      \DOIprefix\doi{10.1145/3067695.3075986}.
    %Type = Article
    \bibitem[{Serigstad et~al.(2018)Serigstad, Eriksen and Breivik}]{Serigstad2018}
    \bibinfo{author}{Serigstad, E.}, \bibinfo{author}{Eriksen, B.O.H.},
      \bibinfo{author}{Breivik, M.}, \bibinfo{year}{2018}.
    \newblock \bibinfo{title}{Hybrid collision avoidance for autonomous surface
      vehicles}.
    \newblock \bibinfo{journal}{IFAC-PapersOnLine} \bibinfo{volume}{51},
      \bibinfo{pages}{1--7}.
    \newblock \URLprefix
      \url{https://linkinghub.elsevier.com/retrieve/pii/S2405896318321499},
      \DOIprefix\doi{10.1016/j.ifacol.2018.09.460}.
    %Type = Article
    \bibitem[{Seta et~al.(2004)Seta, Inoue and USsui}]{Seta2004}
    \bibinfo{author}{Seta, H.}, \bibinfo{author}{Inoue, K.},
      \bibinfo{author}{USsui, H.}, \bibinfo{year}{2004}.
    \newblock \bibinfo{title}{Supporting information for safe ship handling based
      on the concept of potential area of water}.
    \newblock \bibinfo{journal}{The Journal of Japan Institute of Navigation}
      \bibinfo{volume}{111}, \bibinfo{pages}{63--69}.
    \newblock \URLprefix
      \url{https://www.jstage.jst.go.jp/article/jin/111/0/111_KJ00004696699/_article/-char/ja/},
      \DOIprefix\doi{10.9749/jin.111.63}.
    %Type = Article
    \bibitem[{Szlapczynski and Szlapczynska(2017)}]{Szlapczynski2017}
    \bibinfo{author}{Szlapczynski, R.}, \bibinfo{author}{Szlapczynska, J.},
      \bibinfo{year}{2017}.
    \newblock \bibinfo{title}{Review of ship safety domains: Models and
      applications}.
    \newblock \bibinfo{journal}{Ocean Engineering} \bibinfo{volume}{145},
      \bibinfo{pages}{277--289}.
    \newblock \DOIprefix\doi{10.1016/j.oceaneng.2017.09.020}.
    %Type = Book
    \bibitem[{{The Harbor Information Center for the Security of Ship
      Navigation}(2020)}]{nanko2020}
    \bibinfo{author}{{The Harbor Information Center for the Security of Ship
      Navigation}}, \bibinfo{year}{2020}.
    \newblock \bibinfo{title}{Port of Osaka Entrance and Departure Manual}.
    \newblock \bibinfo{edition}{Revised} ed., \bibinfo{publisher}{Osaka ports and
      harbors bureau}.
    %Type = Article
    \bibitem[{Ueno et~al.(2001)Ueno, Nimura, Miyazaki, Fujiwara, Nonaka and
      Yabuki}]{Ueno2001}
    \bibinfo{author}{Ueno, M.}, \bibinfo{author}{Nimura, T.},
      \bibinfo{author}{Miyazaki, H.}, \bibinfo{author}{Fujiwara, T.},
      \bibinfo{author}{Nonaka, K.}, \bibinfo{author}{Yabuki, H.},
      \bibinfo{year}{2001}.
    \newblock \bibinfo{title}{Model experiment and sea trial for investigating
      manoeuvrability of a training ship}.
    \newblock \bibinfo{journal}{Journal of the Society of Naval Architects of
      Japan} \bibinfo{volume}{2001}, \bibinfo{pages}{71--80}.
    \newblock \URLprefix
      \url{http://joi.jlc.jst.go.jp/JST.Journalarchive/jjasnaoe1968/2001.71?from=CrossRef},
      \DOIprefix\doi{10.2534/jjasnaoe1968.2001.71}.
    %Type = Article
    \bibitem[{Ueno et~al.(2017)Ueno, Suzuki and Tsukada}]{Ueno2017}
    \bibinfo{author}{Ueno, M.}, \bibinfo{author}{Suzuki, R.},
      \bibinfo{author}{Tsukada, Y.}, \bibinfo{year}{2017}.
    \newblock \bibinfo{title}{Estimation of stopping ability of full-scale ship
      using free-running model}.
    \newblock \bibinfo{journal}{Ocean Engineering} \bibinfo{volume}{130},
      \bibinfo{pages}{260--273}.
    \newblock \URLprefix
      \url{https://www.sciencedirect.com/science/article/pii/S0029801816305819},
      \DOIprefix\doi{https://doi.org/10.1016/j.oceaneng.2016.12.001}.
    %Type = Article
    \bibitem[{Ueno and Tsukada(2015)}]{Ueno2015}
    \bibinfo{author}{Ueno, M.}, \bibinfo{author}{Tsukada, Y.},
      \bibinfo{year}{2015}.
    \newblock \bibinfo{title}{Rudder effectiveness and speed correction for scale
      model ship testing}.
    \newblock \bibinfo{journal}{Ocean Engineering} \bibinfo{volume}{109},
      \bibinfo{pages}{495--506}.
    \newblock \URLprefix
      \url{https://www.sciencedirect.com/science/article/pii/S0029801815005132},
      \DOIprefix\doi{https://doi.org/10.1016/j.oceaneng.2015.09.041}.
    %Type = Article
    \bibitem[{Vagale et~al.(2021)Vagale, Oucheikh, Bye, Osen, Fossen, Oucheikh,
      Osen and Fossen}]{Vagale2021}
    \bibinfo{author}{Vagale, A.}, \bibinfo{author}{Oucheikh, R.},
      \bibinfo{author}{Bye, R.T.}, \bibinfo{author}{Osen, O.L.},
      \bibinfo{author}{Fossen, T.I.}, \bibinfo{author}{Oucheikh, R.},
      \bibinfo{author}{Osen, O.L.}, \bibinfo{author}{Fossen, T.I.},
      \bibinfo{year}{2021}.
    \newblock \bibinfo{title}{{Path planning and collision avoidance for autonomous
      surface vehicles I: a review}}.
    \newblock \bibinfo{journal}{Journal of Marine Science and Technology}
      \URLprefix \url{https://doi.org/10.1007/s00773-020-00787-6
      http://link.springer.com/10.1007/s00773-020-00787-6},
      \DOIprefix\doi{10.1007/s00773-020-00787-6}.
    %Type = Article
    \bibitem[{Wang and Chin(2016)}]{Wang2016}
    \bibinfo{author}{Wang, Y.}, \bibinfo{author}{Chin, H.C.}, \bibinfo{year}{2016}.
    \newblock \bibinfo{title}{An empirically-calibrated ship domain as a safety
      criterion for navigation in confined waters}.
    \newblock \bibinfo{journal}{Journal of Navigation} \bibinfo{volume}{69},
      \bibinfo{pages}{257--276}.
    \newblock \URLprefix
      \url{https://www.cambridge.org/core/product/identifier/S0373463315000533/type/journal_article},
      \DOIprefix\doi{10.1017/S0373463315000533}.
    %Type = Article
    \bibitem[{Yamanouchi and Fujii(1972)}]{Yamanouchi1972}
    \bibinfo{author}{Yamanouchi, H.}, \bibinfo{author}{Fujii, Y.},
      \bibinfo{year}{1972}.
    \newblock \bibinfo{title}{On the hard core of the effective domain in yokohama
      port}.
    \newblock \bibinfo{journal}{NAVIGATION} \bibinfo{volume}{38},
      \bibinfo{pages}{56--57}.
    \newblock \DOIprefix\doi{10.18949/jinnavib.38.0_56}.
    %Type = Article
    \bibitem[{Yasukawa and Yoshimura(2015)}]{Yasukawa2015}
    \bibinfo{author}{Yasukawa, H.}, \bibinfo{author}{Yoshimura, Y.},
      \bibinfo{year}{2015}.
    \newblock \bibinfo{title}{Introduction of mmg standard method for ship
      maneuvering predictions}.
    \newblock \bibinfo{journal}{Journal of Marine Science and Technology (Japan)}
      \bibinfo{volume}{20}, \bibinfo{pages}{37--52}.
    \newblock \DOIprefix\doi{10.1007/s00773-014-0293-y}.
    %Type = Inproceedings
    \bibitem[{Yoshimura et~al.(2009)Yoshimura, Nakao and Ishibashi}]{Yoshimura2009}
    \bibinfo{author}{Yoshimura, Y.}, \bibinfo{author}{Nakao, I.},
      \bibinfo{author}{Ishibashi, A.}, \bibinfo{year}{2009}.
    \newblock \bibinfo{title}{Unified mathematical model for ocean and harbour
      manoeuvring}, in: \bibinfo{booktitle}{Proceedings of MARSIM2009}, pp.
      \bibinfo{pages}{116--124}.
    \newblock \URLprefix \url{http://hdl.handle.net/2115/42969}.
    
    \end{thebibliography}
\end{document}